\documentclass[11pt,reqno]{amsart}

\usepackage{amsmath,amssymb,amscd,graphicx, amsthm, fullpage,
  enumerate, hyperref} 

\usepackage{mathrsfs}
\usepackage{mathtools}
\usepackage[all]{xy}
\usepackage{txfonts}
\usepackage{tikz}
\usetikzlibrary{matrix,arrows,calc}

\usepackage{color,xcolor}
\definecolor{darkred}{rgb}{1,0,0} 
\definecolor{darkgreen}{rgb}{0,0.8,0}
\definecolor{darkblue}{rgb}{0,0,1}

\hypersetup{colorlinks,
linkcolor=darkblue,
filecolor=darkgreen,
urlcolor=darkred,
citecolor=darkgreen}

\numberwithin{equation}{section}

\theoremstyle{definition}
\newtheorem{theorem}{Theorem}
\numberwithin{theorem}{section}
\newtheorem{proposition}[theorem]{Proposition}
\newtheorem{lemma}[theorem]{Lemma}

\newtheorem{definition}[theorem]{Definition}
\newtheorem{remark}[theorem]{Remark}
\newtheorem{example}[theorem]{Example}
\newtheorem{notation}[theorem]{Notation}

\newcommand{\inv}{^{-1}} 

\newcommand{\toto}{\rightrightarrows}
\newcommand{\id}{\mathrm{id}}
\newcommand{\gen}{\mathsf{gen}}

\newcommand{\Hom}{\mathsf{Hom}}

\newcommand{\Man}{\mathsf{{Man}}}
\newcommand{\Stack}{\mathsf{Stack}}
\newcommand{\GStack}{\mathsf{GeomStack}}
\newcommand{\LieGpd}{\mathsf{LieGpd}}
\newcommand{\Bi}{\mathsf{Bi}}
\newcommand{\sfC}{\mathsf{C}}
\newcommand{\sfD}{\mathsf{D}}

\newcommand{\Lietwo}{\mathsf{Lie2Alg}}
\newcommand{\Vect}{\mathsf{Vect}}
\newcommand{\Vectp}{\mathsf{Vect}'}

\newcommand{\Gpd}{{\sf Gpd}}
\newcommand{\iso}{{\rm{iso}}}
\newcommand{\strict}{{\rm{strict}}}
 
\newcommand{\ESOE}{ \mathscr{E}\mathit{mb}}

\newcommand{\Der}{\mathsf{Der}}

\newcommand{\calA}{{\mathcal A}}
\newcommand{\calB}{{\mathcal B}}

\newcommand{\calX}{{\mathcal X}}

\newcommand{\B}{\mathbb{B}}

\newcommand{\R}{\mathbb{R}}

\newcommand{\X}{\mathbb{X}}

\newcommand{\fg}{\mathfrak g}
\newcommand{\fh}{\mathfrak h}
\newcommand{\fp}{\mathfrak p}
\newcommand{\fs}{\mathfrak s}
\newcommand{\ft}{\mathfrak t}
\newcommand{\fm}{\mathfrak m}
\newcommand{\fn}{\mathfrak n}

\title{Lie 2-algebras of vector fields}

\author{ Daniel Berwick-Evans} 
\email{danbe@illinois.edu}
\author{ Eugene Lerman}
\email{lerman@illinois.edu}

\address{Department of Mathematics, University of Illinois, Urbana,
  IL 61801}

\begin{document}

\begin{abstract}
  We show that the category of vector fields on a geometric stack has
  the structure of a Lie 2-algebra. This proves a conjecture of
  R.~Hepworth.  The construction uses a Lie groupoid that presents the
  geometric stack. We show that the category of vector fields on the
  Lie groupoid is equivalent to the category of vector fields on the
  stack. The category of vector fields on the Lie groupoid has a Lie
  2-algebra structure built from known (ordinary) Lie brackets on
  multiplicative vector fields of Mackenzie and Xu and the global
  sections of the Lie algebroid of the Lie groupoid. After giving a
  precise formulation of Morita invariance of the construction, we
  verify that the Lie 2-algebra structure defined in this way is
  well-defined on the underlying stack.
\end{abstract}

\maketitle

\setcounter{tocdepth}{1}
\tableofcontents

\section{Introduction}  

Vector fields on a Lie groupoid $G$ form a category 
\cite{Hepworth}.  We denote it by $\X(G)$.  The {\em objects} of $\X(G)$
are the multiplicative vector fields of Mackenzie and
Xu~\cite{MackXu}.  These are functors $v:G\to TG$ satisfying
$\pi_G\circ v = \id_G$ where $TG$ denotes the tangent groupoid and
$\pi_G:TG\to G$ is the projection functor.  A {\em morphism}
$\alpha:v \Rightarrow v'$ in this category is a natural transformation
$\alpha$ such that $\pi_G (\alpha(x)) = \id_x$ for every object $x$ of
the groupoid $G$.  The first result of this paper is \\[-6pt]

\noindent {\bf Theorem}~\ref{lem:03.4}.\quad The category of vector
fields $\X(G)$ on a Lie groupoid $G$ is a (strict) Lie 2-algebra.
That is, $\X(G)$ is a category internal to the category of Lie algebras.
\\[-6pt]

\begin{remark}
  When a manifold $M$ is regarded as a discrete Lie groupoid, $\X(M)$
  is the usual Lie algebra of vector fields on $M$ regarded as a
  discrete Lie 2-algebra.
\end{remark}

To every Lie groupoid $G$ there corresponds the stack $\B G$ of
principal $G$-bundles, and Morita equivalent Lie groupoids $G$ and $H$
correspond to isomorphic stacks $\B G$ and $\B H$.  It is natural to
wonder if the Lie 2-algebra $\X(G)$ lives on the stack~$\B G$ in some
appropriate sense. To start, we can ask whether Morita equivalent Lie
groupoids $G$ and $H$ have ``Morita equivalent'' Lie 2-algebras
$\X(G)$ and $\X(H)$. 
More precisely we could ask for the existence of a (2-)functor $\X$
from the bicategory $\Bi$ of Lie groupoids, bibundles and
isomorphisms of bibundles to an appropriate bicategory of Lie
2-algebras that sends Morita equivalences to Morita equivalences.  It
turns out that such a functor is too much to ask for but there is a
functor from a sub-bicategory of $\Bi$.  

The reasons behind this 
fact can already be seen in the case of manifolds.  Recall that there
is no naturally defined functor from the category of manifolds to the
category of Lie algebras that assigns to each
manifold its  Lie algebra  of vector fields.  However if we restrict
ourselves to the category $\Man_\iso$ whose objects are manifolds and
whose morphisms are diffeomorphisms then there is a perfectly well
defined functor with the desired properties.

Getting back to Lie groupoids, recall that $\Bi$ is a localization of
the strict 2-category of Lie groupoids, internal functors, and
internal natural transformations at the class of functors that are 
fully faithful and essentially surjective, i.e., at the essential equivalences. Lie
2-algebras, internal functors and internal natural transformations
form the strict 2-category $\Lietwo_{strict}$, and localizing at the
essential equivalences produces a bicategory~$\Lietwo$ (see
Subsections~\ref{subsec:loc} and \ref{subsec:2} below).  Let
$\Bi_\iso$ be the sub-bicategory of $\Bi$ whose objects are Lie
groupoids, 1-morphisms are (weakly) {\em invertible} bibundles (i.e.,
Morita equivalences) and 2-morphisms are isomorphisms
of bibundles. We recall that a bicategory with invertible 2-morphism
is, by definition, 
a \emph{(2,1)-bicategory}. \\[-6pt]

\noindent {\bf Theorem}~\ref{thm:04.1}.\quad The map $G\mapsto \X(G)$
that assigns to each Lie groupoid its category of vector fields
extends to a functor
\begin{equation*}
\X: \Bi_\iso \to \Lietwo.
\end{equation*}
In particular, if $P:G\to H$ is a Morita
equivalence of Lie groupoids then $\X(P): \X(G)\to
\X(H)$ is a (weakly) invertible 1-morphism of Lie 2-algebras in the
bicategory $\Lietwo$.
 
\begin{remark}
  In the Lie groupoid literature there are two standard constructions
  that associate a Lie algebra to a Lie groupoid: global sections of
  its Lie algebroid and Mackenzie and Xu's multiplicative vector
  fields. The Lie 2-algebra structure on $\X(G)$ is built out of this
  pair of Lie algebras.  At first pass this might seem surprising:
  neither multiplicative vector fields nor sections of Lie algebroids
  are well-behaved under Morita equivalence of Lie groupoids.
  Theorem~\ref{thm:04.1} shows that combining this pair of Lie
  algebras into a Lie 2-algebra gives us an object that {\em is}
  preserved by Morita equivalence.

\end{remark}

How does the existence of the functor in Theorem~\ref{thm:04.1} imply
that the Lie 2-algebra $\X(G)$ ``lives" on the stack $\B G$?  To
answer this, we need to recall the relationship between the bicategory
$\Bi$ and the 2-category $\Stack$ of stacks over the site of smooth
manifolds. The assignment $G\mapsto \B G$ extends to a fully faithful
functor
\[
\B: \Bi \to \Stack.
\]
The essential image of
the functor $\B$ is the 2-category $\GStack$ of geometric stacks.
Restricting the functor $\B$ to the bicategory $\Bi_\iso$ of
groupoids and Morita equivalences gives us an equivalence of
bicategories
\[
\B:\Bi_\iso \to \GStack_\iso,
\]
where $\GStack_\iso$ is the (2,1)-category of geometric stacks,
isomorphisms of stacks (that is, weakly invertible 1-morphisms of
stacks) and 2-morphisms.  By inverting this equivalence $\B$ and
composing it with the functor $\X$ we get a functor
\begin{equation}\label{eq:stack->Lie2}
\GStack_\iso\xrightarrow{\B\inv} \Bi_\iso \xrightarrow{\X} \Lietwo.
\end{equation}
So in particular we get a functorial assignment of a Lie 2-algebra to
every geometric stack, with isomorphic stacks being assigned
``isomorphic'' Lie 2-algebras.

In \cite{Hepworth} Hepworth introduced a category of vector fields
$\Vect(\calA)$ on a stack $\calA$, which is a groupoid. We introduce a groupoid
$\Vectp(\calA)$ equivalent to $\Vect(\calA)$ which is more convenient
for our purposes.  In particular, the assignment
\[
\calA \mapsto \Vectp(\calA)
\]
easily extends to a functor
\[
\Vectp: \GStack_\iso \to \Gpd
\]
where $\Gpd$ is the (2,1)-category of groupoids,  functors and natural
transformations (which are automatically natural isomorphisms).  We  show that the
functor $\Vectp$ is compatible with the functors $\B: \Bi_\iso \to
\GStack_\iso$ and
$\X:\Bi_\iso \to \Lietwo$ in the following sense. \\[-6pt]

\noindent {\bf Theorem}~\ref{thm:5.5.1}.\quad 
The diagram of (2,1)-bicategories and functors
\[
\xy
(-15,8)*+{\GStack_\iso}="1";
(15,8)*+{\Gpd}="2";
(-15, -8)*+{\Bi_\iso}="3";
(15, -8)*+{\Lietwo}="4";
{\ar@{->}^{\B} "3"; "1"};
{\ar@{->}_{u} "4"; "2"};
{\ar@{->}^{\Vectp}"1"; "2"};
{\ar@{->}_{\X}"3"; "4"};
{\ar@{=>} ^{\Upsilon}(2,-2);(-2,2)} ; 
\endxy 
\]
2-commutes.   Here as above $\Gpd$ denotes the (2,1) category of groupoids, functors and natural isomorphisms, and $u:\Lietwo\to \Gpd$ denotes the functor that assigns
to each Lie 2-algebra its underlying groupoid.  The components of the transformation $\Upsilon$ are weakly invertible functors (i.e., equivalences of categories).  In particular for a
geometric stack $\calA$ the category underlying the Lie 2-algebra
$(\X\circ \B\inv)\, (\calA)$
is equivalent to Hepworth's category $\Vect(\calA)$
of vector fields on the stack. \\
\begin{remark}
Consider a geometric stack $\calA$.  The groupoids $\Vect(\calA) $ and $\Vectp(\calA)$ are equivalent.   Let $G = \B \inv (\calA)$. Then the stacks $\B G$ and $\calA$ are isomorphic and consequently the categories $\Vectp(\calA)$ and $\Vectp(\B G)$ are equivalent.  Since the diagram above 2-commutes and the components of $\Upsilon$ are weakly invertible functors, $\Vectp(\calA)$ is equivalent to the groupoid $u (\X (G))$ underlying the Lie 2-algebra $\X(G)$.   Consequently Hepworth's groupoid $\Vect(\calA)$ of vector fields on a stack is equivalent to the groupoid underlying the Lie 2-algebra $\X (\B\inv (\calA))$.   

A different choice $(\B\inv)'$ of the inverse of $\B$ is isomorphic to $\B\inv$.  Consequently the Lie 2-algebras $\X(\B\inv (\calA))$ and $\X((\B\inv)' (\calA)$ are naturally weakly isomorphic and their underlying groupoids are equivalent.
\end{remark}

\subsection*{Related work}
The recent work of Cristian Ortiz and James Waldron~\cite{OW} lies in
a similar circle of ideas.  Recall that an $\mathcal{LA}$-groupoid is
a groupoid object in Lie algebroids.  Given an
$\mathcal{LA}$-groupoid, Ortiz and Waldron introduce its category of
multiplicative sections and show that it carries a natural strict Lie
2-algebra structure in the language of crossed modules of Lie algebras
(which affords an equivalent description of the category of strict Lie
2-algebras).  They show that if two $\mathcal{LA}$-groupoids are
Morita equivalent then the corresponding crossed modules of Lie
algebras are connected by a zig-zag of equivalences.  Furthermore, to
every stack Ortiz and Waldron assign an ordinary Lie algebra and show
that in the case of proper geometric stacks the set underlying this
Lie algebra is in bijective correspondence with isomorphism classes of
vector fields in Hepworth's definition.

\subsection*{Outline of the paper}
In Section~\ref{sec:2} we review some of the background material used
in the paper.  In particular we recall the strict 2-category $\LieGpd$
of Lie groupoids, smooth functors and natural transformations.  We
then briefly discuss the bicategory $\Bi$ of Lie groupoids, bibundles
and isomorphisms of bibundles and the functor $\langle \,\rangle
:\LieGpd\to \Bi$ that localizes the strict 2-category $\LieGpd$ at the
class of the essential equivalences.  We then discuss the
localizations of bicategories in general and recall a criterion due to
Pronk for a functor between bicategories to be a localization.  We
then review 2-vector spaces, strict Lie 2-algebras and crossed modules
of Lie 2-algebras.  We localize the strict 2-category
$\Lietwo_\strict$ of Lie 2-algebras, internal functors and natural
transformations at essential equivalences and obtain a bicategory
$\Lietwo$. Under the correspondence between Lie 2-algebras and crossed
modules $\Lietwo$ corresponds to Noohi's bicategory of crossed-modules
and butterflies~\cite{Noohi}.  We finish the section by discussing the
extension of the tangent functor $T$ on the category of manifolds to
tangent functors on the 2-category $\LieGpd$ and the bicategory $\Bi$,
respectively.

In Section~\ref{sec:3} we prove Theorem~\ref{lem:03.4}: the category
of multiplicative vector fields on a Lie groupoid underlies a strict
Lie 2-algebra.  In Section~\ref{sec:4} we prove that the assignment
$G\mapsto \X(G)$ of the category of vector fields to a Lie groupoid
extends to a functor $\X: \Bi_\iso \to \Lietwo$ from the bicategory
$\Bi_\iso$ of Lie groupoids, invertible bibundles and isomorphisms of
bibundles to the bicategory $\Lietwo$ of Lie 2-algebras.  Hence, in
particular, if $P:G\to H$ is a Morita equivalence of Lie groupoids
then $\X(P): \X(G)\to \X(H)$ is a (weakly) invertible 1-morphism of
Lie 2-algebras in the bicategory $\Lietwo$.  Along the way we
introduce the category $\X_\gen(G)$ of {\em generalized } vector
fields on a Lie groupoid $G$.  It is modelled on Hepworth's category
of vector fields on a stack. The objects of $\X_\gen(G)$ are pairs
$(P, \alpha_P)$ where $P:G\to TG$ is a bibundle and $\alpha_P:\tilde{
\pi}_G \circ P\Rightarrow \langle \id_G\rangle$ is an
isomorphism of bibundles.  Here and below $\tilde{ \pi}_G$ is
the bibundle corresponding to the projection functor $\pi_G:TG\to G$
(see Lemma~\ref{lem:2.26'})
and $\langle \id_G \rangle$ is the identity bibundle on the Lie
groupoid $G$.

In Section~\ref{sec:5} we discuss Hepworth's category of vector fields
$\Vect(\calA)$ on a stack $\calA$ and construct an equivalent category
$\Vectp(\calA)$. In Section~\ref{sec:5.5} we promote the assignment
$\calA \mapsto \Vectp(\calA)$ of a category of vector fields on a
geometric stack to a functor $\Vectp:\GStack_\iso \to \Gpd$ from the
2-category of geometric stacks and isomorphisms to the  2-category $\Gpd$
of groupoids and prove Theorem~\ref{thm:5.5.1}.  In
Section~\ref{sec:6} we prove Theorem~\ref{thm:04.3}: for any Lie
groupoid $G$ the ``inclusion'' functor
\[
\imath_G: \X(G)\hookrightarrow \X_\gen(G), \qquad v\mapsto (\langle v\rangle,
\alpha_{\langle v\rangle }:\langle \pi_G\rangle \circ \langle v
\rangle \Rightarrow \langle id_G \rangle).
\]
of the category of multiplicative vector fields into the category of
generalized vector fields is fully faithful and essentially
surjective.  This generalizes a result of Hepworth for proper Lie
groupoids.

\subsection*{Acknowledgments} 
We thank Henrique Bursztyn for many helpful discussions.  In
particular this paper has partially originated from conversations of
one of us (E.L.) with Henrique at Poisson 2014. We thank James Waldron
for making us exercise more care with the equivalences of Lie
2-algebras.

We thank the anonymous referees for their helpful comments.

\section{Background and notation}\label{sec:backgr}\label{sec:2}

We assume that the reader is familiar with ordinary categories, strict
2-categories and bicategories (also known as weak 2-categories).  We
mostly 
work  with (2,1)-bicategories, that is with bicategories
whose 2-morphisms are invertible. Standard references for bicategories
are \cite{Benabou} and \cite{Borceux}.  For the reader's convenience
the definitions of a bicategory, (pseudo-)functors, (pseudo-natural)
transformations, modifications and functor bicategories  are summarized in Appendix~\ref{sec:app}.   We assume familiarity with Lie
groupoids.  Standard references are \cite{MacK} and \cite{MM}.  We
also assume that the reader is comfortable with stacks over the site
of manifolds.  This said, sections \ref{sec:3}, \ref{sec:4} and
\ref{sec:6} do not use stacks. 
  While there is no
textbook covering stacks over manifolds, a number of references exist:
\cite{BehXu}, \cite{Blo}, \cite{L}, \cite{Metzler}, \cite{SP} (this list is not
exhaustive).

Given a category $\sfC$ we denote its collection of objects by
$\sfC_0$ and the collection of arrows/morphisms by
$\sfC_1$.\footnote{We use the words ``arrow,'' 
  ``morphism'' and ``1-cell'' interchangeably.}  We usually denote the source and
target maps of $\sfC$ by $s$ and $t$, respectively.  We write
\[
\sfC = \{\sfC_1\toto \sfC_0\}
\]
and suppress the other structure maps of the category $\sfC$.  We
denote the unit map by $1$.  Thus the map $1:\sfC_0\to \sfC_1$ assigns
the identity arrow $1_x$ to each object $x$ of the category $\sfC$.
The composition/multiplication in the category $\sfC$ is defined on
the collection $\sfC_2$ of pairs of composable arrows.  Our convention
is that
\[
\sfC_2:= \{(\gamma_2, \gamma_1) \in \sfC_1 \times \sfC _1 \mid
s (\gamma_2) = t (\gamma_1)\} =: \sfC_1\times_{s, \sfC _0, t} \sfC_1.
\]
We denote the composition in the category $\sfC$ by $m$:
\[
m:\sfC_1\times _{s, \sfC _0,t } \sfC_1  \to \sfC_1, \qquad
(\gamma_2, \gamma_1)\mapsto m(\gamma_2, \gamma_1) \equiv \gamma_2\gamma_1.
\]
In particular, we write the composition from right to left:
$\gamma_2\gamma_1$ means $\gamma_1$ followed by $\gamma_2$.
If the category $\sfC$ is a groupoid we denote the inversion map by $i$:
\[
i:\sfC_1\to \sfC_1,\qquad i(\gamma) := \gamma\inv.
\]

\subsection{Bicategories of Lie groupoids}

We start with recalling the standard  notion of a bibundle between two Lie groupoids.
Bibundles are also known as generalized morphisms and as 
Hilsum-Skandalis maps.

\begin{definition}
A {\sf bibundle} $P:G\to H$ from a Lie groupoid $G$ to a Lie groupoid
$H$ is a manifold $P$ with two maps  $a^L_P$ and
$a^R_P$:
\[
\xy
(-19,8)*+{G_1}="1";
(-19,-10)*+{G_0}="2";
(19, 8)*+{H_1}="3";
(19, -10)*+{H_0}="4";
(0,5)*+{P} = "5";
{\ar@{->} (-20,6); (-20, -7)};
{\ar@{->} (-18,6); (-18, -7)};
{\ar@{->} (20,6); (20, -7)};
{\ar@{->} (18,6); (18, -7)};
{\ar@{->}_{a^L_P}  "5";"2"};
{\ar@{->}^{a^R_P}  "5";"4"};
\endxy\qquad 
\]
along with a left action of $G$ and a right action of $H$ on $P$:
\[
G_1\times _{G_0}P\to P\quad (g,p)\mapsto g\cdot p , \qquad P\times
_{H_0}H_1 \to P\quad (p,h)\mapsto p\cdot h.
\]
We refere to $a_P^L$ as the {\sf left anchor} and to $a_P^R$ as the {\sf right
anchor}.  We further require that
\begin{enumerate}
\item the actions of $G$ and $H$ commute: for all $(g,p)\in
  G_1\times_{G_0} P$, $(p, h)\in P\times _{H_0} H_1$
  \[
g\cdot (p\cdot h) = (g \cdot p)\cdot h;
\]
\item the map $a_L^P:P\to G_0$ is a surjective submersion and is
  $H$-invariant: $a_L^P (p\cdot h) = a_L^P (p)$ for all $(p, h)\in
  P\times _{H_0} H_1$;
  \item the map $a^P_R$ is $G$-invariant: $a^P_R (g\cdot p) = a^P_R
    (p)$ for all $(g,p)\in
  G_1\times_{G_0} P$;
\item  the map
  \[
P\times_{a^R_P,H_0,t} H_1\to P\times_{a^L_P,G_0,a^L_P}P \qquad
(p,h)\mapsto (p, p\cdot h)
\]
is a diffeomorphism.
  \end{enumerate}
\end{definition}
\begin{remark}
The conditions on the bibundle $P:G\to H$ guarantee that the action of
$H$ on $P$ makes the submersion $a^P_L:P\to
G_0$ into  a principal $H$-bundle.
  \end{remark}

\begin{definition}
  An {\sf isomorphism} of two bibundles $P,P':G\to H$ is a
  diffeomorphism $\alpha: P\to P'$ which is $G$- and $H$-equivariant.
  \end{definition}
Bibundles can be composed.   Given two bibundles 
$P:G\to H$ and $Q:H\to K$  their composite  $Q\circ P$ is defined to
be the quotient of the fiber product $P_{a^P_R, H_0, a^Q_L} Q$ by the
action of $H$:
\begin{equation}
Q\circ P := (P\times _{a^P_R, H_0, a^Q_L} Q)/H \label{eqn:bibuncomp}.
\end{equation}
The composition of bibundles is not associative on the nose: given 3
composable bibundles $P, Q$ and $R$ the composites $(R\circ Q)\circ P$
and $R\circ (Q\circ P)$ are only isomorphic.  Therefore Lie groupoids
and bibundles do not form a category.
One can show that  Lie groupoids, bibundles and isomorphisms of bibundles
form a bicategory which we denote by $\Bi$.  See \cite{Blo} or
\cite{L}.

\begin{notation}
An isomorphism of bibundles $\alpha: P\to Q$ is a smooth map and
2-arrow (2-cell) in the bicategory $\Bi$ described above.  For this
reason we may sometimes write $\alpha: P\Rightarrow Q$ for an
isomorphism of bibundles.  
  \end{notation}

There is also the strict 2-category $\LieGpd$ of Lie groupoids, (smooth)
  functors and (smooth) natural transformations (which are
  automatically natural isomorphisms).   Note that  both $\Bi$ and
  $\LieGpd$ are (2,1) bicategories.  
    


\begin{notation}
  In the bicategories $\LieGpd$ and $\Bi$, we write the horizontal
  composition of 2-arrows as $\star$. Given a 1-morphism $f$ and a
  2-morphism $\alpha$ we abuse notation by writing $f\star\alpha$ for
  the horizontal composition (whiskering) $1_f \star \alpha$ where
  $1_f$ is the identity 2-arrow on the 1-morphism $f$.  The vertical
  composition of 2-morphisms is denoted by $\circ$.  When convenient,
  we also use arrow notation to denote morphisms in groupoids with
  specified source or target, e.g., $x\xleftarrow{g} y$ for a morphism
  $g$ with target $x$ and source~$y$.
\end{notation}

\begin{remark} \label{rmrk:02.2}
There is a functor (for example, see \cite{L})
\begin{equation}
U: \LieGpd \to \Bi 
\end{equation}
that is the identity on objects.  On 
1-morphisms $U$ sends a functor $f:G\to H$ to the bibundle
\begin{equation} \label{eq:02.3} 
\langle f \rangle := G_0\times_{f,
    H_0, t} H_1 := \{(x,\gamma) \mid f(x) = t(\gamma)\} = \{ (x, f(x)
  \xleftarrow{\gamma})\mid x\in G_0, \gamma \in H_1\}
\end{equation}
whose left and right anchor maps are given respectively by 
\[
a^L_{\langle f\rangle} (x, \gamma) = x, \qquad a^R_{\langle f \rangle}(x,\gamma) = s(\gamma). 
\]
Here, as before, $s:H_1\to H_0$ is the source map.  The left action of
the groupoid $G$ on the manifold $\langle f\rangle$ is given by
\[
(g, (x, \gamma)) \mapsto (t(g), f(g) \gamma).
\]
The right action of the groupoid $H$ on  $\langle f\rangle$ is given by
\[
((x, \gamma), \nu ) \mapsto (x,  \gamma \nu).
\]
Note that $a^L_P: \langle f \rangle \to G_0$ has a canonical section
\[
x \mapsto (x, 1_{f(x)}).
\]
Given a pair of functors $f,k:G\to H$ and a natural
isomorphism $\alpha:f\Rightarrow k$, we get an isomorphism of
bibundles
\[
\langle \alpha \rangle : \langle f\rangle \Rightarrow \langle k\rangle.
\]
The isomorphism $\langle \alpha \rangle$ is defined by
\[
\langle \alpha \rangle (x, f(x) \xleftarrow{\gamma}) = (x, k(x)
\xleftarrow{\alpha (x)\gamma}).
\]
It is not hard to check that the map $\langle \alpha \rangle$ defined above is
smooth, commutes with the left and right anchor maps and is 
equivariant with respect to the actions of $G$ and $H$.

The functor $U= \langle \,\, \rangle$ takes vertical
and horizontal composition of natural transformations to the
composition of isomorphisms of bibundles and horizontal composition of
isomorphisms, respectively.   
\end{remark}
\begin{remark}
  By construction of the functor $U= \langle \,\,\rangle$ the total space
  of the bibundle $\langle id_G\rangle$ corresponding to the identity
  functor $id_G:G\to G$ on a Lie groupoid $G$ is the fiber product
  $G_0\times _{G_0}G_1$.  This fiber product is diffeomorphic
  to~$G_1$.  We therefore define the manifold $G_1$ together with the
  actions of $G$ by left and right multiplication to be the identity
  bibundle for a Lie groupoid $G$.
\end{remark}%

The functor $U$ is far from being an equivalence of
2-categories.  The issue is that for almost all groupoids $G$ and $H$
the functor
\begin{equation} \label{eq:not-ess-surj}
U: \Hom_\LieGpd (G, H) \to \Hom_\Bi( G, H)
\end{equation}
fails to be essentially surjective. The failure of essential
surjectivity follows from the well-known fact:
\begin{lemma}\label{lem:2.5n}
  A bibundle $P:G\to H$ is isomorphic to a bibundle
  $\langle f \rangle$ for some functor $f:G\to H$ if and only if the
  left anchor $a_P^L:P\to G_0$ has a section.
\end{lemma}
\begin{proof}
  For a functor $f:G\to H$, the left anchor map $a^L_P: \langle f\rangle
  = G_0\times_{H_0}H_1 \to G_0$ has a canonical global section $x\mapsto (x,
  1_{f(x)})$.

  Conversely suppose $a^L_P: P\to G_0$ has a global section $\sigma$.  We
  define the corresponding functor $f_\sigma:G\to H$ on objects by
\[
f_\sigma(x) := (a^R_P \circ \sigma)(x).
\]
Since $a^L_P:P\to G_0$ is a principal $H$-bundle for any arrow
$y\xleftarrow{\gamma} x \in G_1$ there is a unique arrow $\tau$ in $H_1$, 
which depends smoothly on $\gamma$, so that
\[
\gamma \cdot \sigma(x) = \sigma(y) \cdot \tau.
\]
We set
\[
f_\sigma (\gamma):=\tau.
\]
It is easy to check that $f$ is indeed a morphism of Lie groupoids.
\end{proof}

\begin{remark}
 Recall that a functor $f:G\to
H$ between two Lie groupoids is an {\sf essential equivalence} if the
right anchor $a^R: \langle f\rangle\to  H_0$ is a surjective submersion
and the action of $G$ on $\langle f \rangle$ is principal.  The functor
$U =:\LieGpd \to \Bi$ is a localization of the 2-category
$\LieGpd$ at the class of all  essential equivalences; see Example~\ref{ex:2.12} below.
\end{remark}
In contrast to failure of the functor  $U$ to be surjective on
1-morphisms, for 2-morphisms the following result holds.  The result
must be known but we are not aware  of a  reference.
\begin{theorem}[Folklore] \label{thm:folk1}
For any pair of functors $f,k:G\to
  H$ of Lie groupoids the map
\[
U: \Hom_\LieGpd(f,k) \to 
\Hom_\Bi(\langle f\rangle, \langle k\rangle),
\qquad \alpha \mapsto \langle \alpha \rangle 
\]
is a bijection.
\end{theorem} 

\begin{proof}[Sketch of proof] Let $\delta: \langle f \rangle \to
  \langle k\rangle$ be an isomorphism of bibundles.  The left anchor
  $a^L_{\langle f\rangle }:\langle f\rangle \to G_0$ has a natural
  section $\sigma_f$.  It is defined by
\[
\sigma_f (x) = (x, f(x)\xleftarrow{1_{f(x)}}).
\]
Similarly we have a natural section $\sigma_k:\langle k\rangle \to
G_0$ of the left anchor $a^L_{\langle k\rangle }:\langle k\rangle \to
G_0$.  Since $a^L_{\langle k\rangle }:\langle k\rangle \to G_0$ is a
principal $H$ bundle, for any $x\in G_0$ there is a unique arrow
$\bar{\delta}(x)\in H_1$ so that
\[
\delta (\sigma_f(x))= \sigma_k (x)\cdot \bar{\delta}(x)
\]
for all $x\in G_0$.  By equivariance of $\delta$, the map
\[
\bar{\delta}:G_0 \to H_1\qquad x\mapsto \bar{\delta}(x)
\]
is a natural isomorphism from $f$ to $k$.
\end{proof}
\subsection{Localizations of bicategories}\label{subsec:loc}
Appendix~\ref{sec:app} recalls definitions of functors,
transformations and modifications in the context bicategories.  We
also recall   functor bicategories.
\begin{definition}
Given a bicategory $\sfC$ and a class of 1-morphisms $W$ in $\sfC$ we
define a  {\sf
  localization}  of $\sfC$ at the class $W$ (if it exists) to be a pair
$(\sfC[W\inv], U:\sfC \to \sfC[W\inv])$ where 
$\sfC[W\inv]$ is a bicategory and $U:\sfC \to \sfC[W\inv]$ is a
functor 
satisfying the following universal property:  For any bicategory $\sfD$
the precomposition with $U$ induces an equivalence of bicategories
\[
\Hom (\sfC[W\inv], \sfD) \xrightarrow{-\circ U} \Hom _W (\sfC, \sfD)
\]
where $\Hom _W (\sfC, \sfD)$ denotes the bicategory of functors sending
elements of $W$ to weakly invertible 1-morphisms in $\sfD$.
\end{definition}
In
particular given any functor $F:\sfC \to \sfD$ mapping elements of $W$
to invertible morphisms of $\sfD$ there exists a functor
$\tilde{F}:\sfC[W\inv] \to \sfD$ and a natural isomorphism
\[
F \Rightarrow \tilde{F}\circ U.
\]
Following a common abuse of notation, we often denote the localizations of $\sfC$
at $W$ simply as
~$\sfC[W\inv]$ (and omit the localization functor $U$).

The localization $\sfC[W\inv]$ is defined up to equivalence of
bicategories, so it will be convenient to refer to any functor $F:\sfC
\to \sfC'$ between bicategories as a localization of $\sfC$ at the
class $W$ if it has the same universal property as the localization functor $U:\sfC \to
\sfC[W\inv]$.  Namely we ask that for any bicategory $\sfD$ the
precomposition with $F$ induces an equivalence of bicategories
\[
\Hom (\sfC', \sfD) \xrightarrow{-\circ F} \Hom _W (\sfC, \sfD).
\]
Pronk \cite{Pronk} gives a criterion for a functor $F:\sfC \to
\sfC'$ to be ``the'' localization of $\sfC$ at the class $W$:
\begin{proposition} \label{prop:Pronk}
(\cite[Proposition~24]{Pronk}) A functor $F:\sfC
  \to \sfC'$ between bicategories is a localization of $\sfC$ at
  the class $W$ if
 \begin{enumerate}
 \item $F$ sends the elements of $W$ to (weakly) invertible 1-morphisms
   in $\sfC'$;
\item $F$ is essentially surjective on objects;
\item for every 1-morphism $f$ in $\sfC'$ there are 1-morphisms $w$ in
  $W$ and $g$ in $\sfC$ with a 2-morphism $F(g) \Rightarrow f\circ F(w)$;
\item $F$ is fully faithful on 2-morphisms.
\end{enumerate}
\end{proposition}

\begin{example} \label{ex:2.12} The functor $U:= \langle \,\rangle:\LieGpd \to \Bi$ is the
  localization of the 2-categories of Lie groupoids, functors and
  natural transformations at the class of essential equivalences.
  Indeed the functor is surjective on objects and sends essential
  equivalences to invertible bibundles.  Finally, for any bibundle $P$
  a choice of local sections of the left anchor leads to a
  factorization $P \Rightarrow \langle g \rangle \circ \langle w
  \rangle\inv$ where $g$ is a functor and $w$ is an essential
  equivalence.
\end{example}

Localizations of bicategories will come up several times in this
paper.  For example, in the next subsection we will discuss the localization of the
strict category $\Lietwo_{strict}$ of Lie 2-algebras at essential
equivalences.  In Section~\ref{sec:4} we will need the fact that the
bicategory $\Bi_\iso$ of Lie groupoids, weakly invertible bibundles
and isomorphisms of bibundles is a localization of a certain 2-category
of embeddings of Lie groupoids.

\subsection{2-vector spaces and Lie 2-algebras} \label{subsec:2}
\begin{definition}
A {\sf 2-vector space} (in the sense of Baez and Crans \cite{BC}) is a
category $V$ internal to the category of vector spaces.  Hence
$V=\{V_1\toto V_0\}$ where $V_0$ a vector space of objects, $V_1$ a
vector space of morphisms, and all the structure maps (source, target,
unit, and composition) are linear.  All 2-vector spaces in this paper
are defined over~$\R$.
\end{definition}

There is a 2-category $2\mathsf{Vect}$ whose objects are 2-vector spaces, 
1-morphisms are (linear) functors and 2-morphisms are
(linear) natural transformations.  There is a forgetful functor 
\begin{equation}
2\mathsf{Vect} \to \mathsf{Cat}
\end{equation}
from the 2-category of 2-vector spaces to the 2-category $\mathsf{Cat}$
of categories that forgets the linear structure.

\begin{remark} \label{rmrk:2term-2vect_corr}
  There is an equivalence of categories of 2-vector spaces and of
  2-term chain complexes of vector spaces. See, for example,
  \cite{BC}.  A similar result characterizing Picard categories was
  obtained much earlier by Deligne \cite{SGA4}.   We remind the reader
  of how this equivalence is defined on objects.  Given a 2-term
  complex $\partial: U\to W$ there is an action of the abelian group
  $U$ on $W$ given by
\begin{equation}\label{eq:comp_in_T}
u\cdot w:= \partial(u) + w
\end{equation}
for all $u\in U, w\in W$. The corresponding action groupoid $\{U\times
W\toto W\}$ is a 2-vector space.  

The converse is true as well: any 2-vector space $V=\{V_1\toto V_0\}$
is isomorphic to an action groupoid defined by the 2-term complex
$\partial = t|_{\ker s} :\ker s \to V_0$. Here as before $s, t:V_1\to
V_0$ are the source and target map of the category $V$; see \cite{BC}
for a proof.   In particular a category underlying a 2-vector
space is a groupoid.
\end{remark}

Next we recall the definition of a strict Lie 2-algebra \cite{BC}.
\begin{definition}
  A strict {\sf Lie 2-algebra} is a  category internal to the
  category of Lie algebras (over the reals): the space of objects and
  morphisms of a Lie 2-algebra are ordinary Lie algebras and all the
  structure maps are maps of Lie algebras.
\end{definition}

\begin{notation}
  Categories internal to Lie algebras, internal functors and internal
  natural transformations form a strict 2-category which we denote by
  $\Lietwo_{strict}$.
\end{notation}

\begin{definition} \label{def:cross} (see, for example,
  \cite[Definition~15]{FW}) A {\sf crossed module} of Lie algebras
  consists of a Lie algebra homomorphism $\partial : \fm \to \fn$
  together with a Lie algebra homomorphism
\[
D: \fn \to \Der(\fm)%
\]
from $\fn$ to the Lie algebra $\Der(\fm)$ of derivations of $\fm$
so that for all $m,m' \in \fm$, $n\in \fn$
\begin{itemize}
\item[(i)\,] $\partial (D(n) m) = [n, \partial (m)]$ and 
\item[(ii)\,] $D(\partial (m))m' = [m,m']$.
\end{itemize}
\end{definition}

A crossed module of Lie algebras determines a Lie 2-algebra: see, for
example, the proof of Theorem~3 in \cite{FW}.  The converse is true as
well: any Lie 2-algebra canonically defines a crossed module of Lie
2-algebras.  In fact more is true: crossed modules form a strict 
2-category, and the 2-categories of Lie 2-algebras and of crossed
modules are equivalent (see \cite[Theorem~3]{FW} cited above).  We
won't need the full strength of this theorem in the present paper.  We
do, however, need the following result:

\begin{lemma}\label{lem:Lie2_from_cross}
  Let $V=\{V_1\toto V_0\}$ be a 2-vector space.  Suppose the
  corresponding 2-term complex $\partial = t|_{\ker s}:\ker s \to V_0$
  is part of the data of a Lie algebra crossed module.  That is,
  suppose that $V_0$, $\ker s$ are Lie algebras, $\partial$ is a Lie
  algebra map, and that there is an action $D:V_0 \to \Der (\ker s)$
  of $V_0$ on $\ker s$ by derivations making $(\partial:\ker s\to V_0,
  D:V_0 \to \Der (\ker s))$ into a crossed module of Lie
  algebras. Then $V$ is a Lie 2-algebra.
\end{lemma}

\begin{proof}[Sketch of proof] 
  Since $1\circ s = id_{V_0}$, $V_1 = \ker s \oplus V_0$.  We
  define a bracket on $\ker s \oplus V_0$ by
\begin{equation}\label{eq:v_1brack}
[(x_1, y_1), (x_2, y_2)] := ([x_1, x_2]+ D(y_1) x_2 - D(y_2)x_1, [y_1,y_2])
\end{equation}
for all $(x_1, y_1), (x_2, y_2) \in \ker s \oplus V_0$. That is, we
define the Lie algebra $V_1$ to be the semi-direct product of $V_0$ and
$\ker s$.  Checking that source, target and unit maps of $V$ are Lie
algebra maps is easy.  To check that the composition $m:V_1 \times
_{V_0}V_1\to V_1$ in the category $V$ is a Lie algebra map we observe
 that $m$ is given by
\begin{equation}\label{eq:m}
m((x_1, y_1),(x_2, y_2))=(x_1+ x_2, y_2)
\end{equation}
for all $(x_1, y_1), (x_2, y_2) \in \ker s \oplus V_0$ with $y_1 =
t(x_2, y_2) = \partial x_2 + y_2$.  This fact is not completely
obvious.  It lies in the heart of the correspondence between 2-vector
spaces and 2-term chain complexes.  See
Remark~\ref{rmrk:2term-2vect_corr} and \cite{BC}.  A computation now
shows that the map $m$ is a Lie algebra map.
\end{proof}

There is a problem with the 2-category $\Lietwo_{strict}$ of Lie
2-algebras.  Namely, suppose $f:\fg\to \fh$ is a morphism of Lie
2-algebras which is fully faithful and essentially surjective, that
is, an essential equivalence.  Then $f$ has a weak inverse (as a
functor), but there is no reason for that inverse to be a morphism of
Lie 2-algebras.  In fact it is easy to come up with examples where
such morphism of Lie 2-algebras does not exist.  Here is one.  The
$2n+1$ dimensional Heisenberg Lie algebra $\mathfrak{h}$ is a central
extension of a $2n$ dimensional abelian Lie algebra $\mathfrak{a}$ by
the reals. Consequently we have a map $\varphi$ of 2-term complexes of
Lie algebras $(\R \to \mathfrak{h}) \to (0\to \mathfrak{a})$, but the
map $\phi_0: \mathfrak{h}\to \mathfrak{a} $ has no Lie algebra
sections.  In fact $\varphi$ is a map of crossed modules of Lie
algebras.  The morphism $\varphi$ of crossed modules corresponds to an
essential equivalence of Lie 2-algebras for which there is no inverse
map in $\Lietwo_{strict}$.

Fortunately the problem has a universal solution: we localize the
2-category $\Lietwo_{strict}$ at the class of essential equivalences and
obtain a bicategory $\Lietwo$ (see \cite{Pronk} and
Subsection~\ref{subsec:loc} above).  This localization has a simple
and explicit description: we define a morphism between Lie 2-algebras as
a ``bibundle internal to the category of Lie
algebras.''  
 Here are the details.  

\begin{definition} \label{def:lie-bi}
A {\sf bibundle} $\fp: \fg
  \to \fh$ from a Lie 2-algebra $\fg$ to a Lie 2-algebra
  $\fh$ is a Lie algebra $\fp$ with {\sf left} and {\sf right anchor
    maps} $a_\fp^L$ and $a_\fp^R$ (which are maps of Lie algebras),
\[
\xy
(-19,8)*+{\fg_1}="1";
(-19,-10)*+{\fg_0}="2";
(19, 8)*+{\fh_1}="3";
(19, -10)*+{\fh_0}="4";
(0,5)*+{\fp} = "5";
{\ar@{->} (-20,6); (-20, -7)};
{\ar@{->} (-18,6); (-18, -7)};
{\ar@{->} (20,6); (20, -7)};
{\ar@{->} (18,6); (18, -7)};
{\ar@{->}_{a^L_\fp}  "5";"2"};
{\ar@{->}^{a^R_\fp}  "5";"4"};
\endxy
\]
along with a left action of the groupoid $\fg$ and right action of the
groupoid $\fh$
\[
\fg_1\times_{s,\fg_0,a_\fp^L}\fp\to \fp \quad (g,p)\mapsto g\cdot p, \quad \quad 
\fp\times_{a_\fp^R,\fh_0,t}\fh_1\to \fp\quad (p,h)\mapsto p\cdot h.
\]
We require that the actions are maps of Lie algebras and commute with
each other.  We require that $a_\fp^L$ is surjective. Finally, we
require that the map
\[
\fp\times_{a_\fp^R,\fh_0,t}\fh_1
\to \fp\times_{a_\fp^L,\fg_0,a_\fp^L} \fp\quad \quad (p,h)\mapsto (p,p\cdot h)
\]
is an isomorphism of Lie algebras.  Thus in particular  we require that $a^L_\fp :\fp
\to \fg_0$ is a principal 
$\fh$-bundle.
\end{definition}

\begin{remark}\label{rmrk:02.13}
  The composition of bibundles between Lie 2-algebras is defined
  in the same way as in the case of bibundles between Lie groupoids:
  it is the quotient of the appropriate fiber product.
  We will omit a proof that Lie 2-algebras, bibundles of Lie algebras
  and isomorphisms of bibundles form a bicategory. We denote this bicategory by
  $\Lietwo$.  We note that biprincipal bibundles are weakly invertible
  in this bicategory.

  As in the case of Lie groupoids there is a functor $\langle \,\,
  \rangle :\Lietwo_{strict} \to \Lietwo$. It sends a strict map
  $f:\fg\to \fh$ of Lie 2-algebras to the bibundle
\[
\langle f \rangle := \fg_0\times_{f_0, \fh_0, t} \fh_1 := \{(x,\gamma) \mid
f_0(x) = t(\gamma)\},
\]
whose left and right anchor maps are given respectively by 
\[
a^L (x, \gamma) = x, \qquad a^R(x,\gamma) = s(\gamma). 
\]
The left action of $\fg$ on $\langle f\rangle$ is 
\[
(g, (x, \gamma)) \mapsto (t(g), f(g) \gamma),
\]
and the right action of $\fh$ on  $\langle f\rangle$ is 
\[
((x, \gamma), \nu ) \mapsto (x,  \gamma \nu).
\]
\end{remark}

In order to check that the functor $\langle \quad \rangle
:\Lietwo_{strict} \to \Lietwo$ is the localization of the 2-category
$\Lietwo_{strict}$ at the class of essential equivalences we need

\begin{lemma}\label{lem:weak-inv}
  Suppose $f:\fg \to \fh$ is a strict map of Lie 2-algebras whose
  underlying functor is fully faithful and essentially surjective.
  Then the bibundle of Lie 2-algebras
\[
\langle f \rangle :\fg \to \fh
\]
is weakly invertible.
\end{lemma}

\begin{proof}
  It is enough to show that $a^R_{\langle f \rangle}: \langle f
  \rangle \to \fh$ is an
  $\fg$-principal bundle.  That is, it is enough to show that $a^R_{\langle f \rangle} $ is
  surjective and that the map
\[
\varphi: \fp\times_{a_\fp^R,\fh_0,t}\fh_1%
\to \fp\times_{a_\fp^L,\fg_0,a_\fp^L} \fp\quad \quad \varphi(p,h):= (p,p\cdot h)
\]
is an isomorphism of Lie algebras. Since $a^R_{\langle f \rangle}(x,\gamma) = s(\gamma)$
the surjectivity of $a^R_{\langle f \rangle}$ is equivalent to the essential surjectivity
of the functor $f$.  The fullness of $f$ translates into $\varphi$
being onto and faithfulness of $f$ translates into $\varphi$ being
1-1.
\end{proof}
We now apply Proposition~\ref{prop:Pronk} to conclude that $\langle
\quad \rangle :\Lietwo_{strict} \to \Lietwo$ is the localization of
$\Lietwo_{strict} $ at the class of essential equivalences.  See also
Theorem~\ref{thm:loc} below for a similar argument.

\begin{remark}
  A reader familiar with Noohi's butterflies (see \cite{Noohi} and
  reference therein) should not have much trouble showing that the
  bicategory $\Lietwo$ of Lie 2-algebras defined above is equivalent to the
  bicategory of crossed modules of Lie algebras, butterflies and
  isomorphisms of butterflies.

  Alternatively, this equivalence can be
  seen as an equivalence of localizations of two  equivalent
  2-categories. Indeed, as recalled above the strict 2-category of Lie 2-algebras is 
  equivalent to the 2-category of crossed modules of Lie algebras. 
  Noohi's butterflies localize crossed modules at the class of 1-morphisms that correspond
  precisely to the class of essential equivalences of Lie 2-algebras.
\end{remark}

\subsection{Tangent functors} 
\label{subsec:tangent}

Recall that we have a functor $T$ from the category $\Man$ of
$C^\infty$ manifolds to itself. 
To each map $f:M\to N$ between manifolds the functor assigns the differential
$Tf:TM\to TN$ between their tangent bundles. The fact that $T$ preserves the
composition of maps is the chain rule.  Moreover for any map $f:M\to
N$  the diagram
\[
\xy
(-10,10)*+{TM}="1";
(10,10)*+{TN}="2";
(-10, -6)*+{M}="3";
(10, -6)*+{N}="4";
{\ar@{->}^{Tf} "1"; "2"};
{\ar@{->}_{\pi_M} "1"; "3"};
{\ar@{->}^{\pi_N} "2"; "4"};
{\ar@{->}_f "3"; "4"};
\endxy
\]
commutes.  Hence we have a natural transformation $\pi:T\Rightarrow
\id_{\Man}$ whose components are the projections $\pi_M:TM\to M$.

\begin{remark}\label{rmrk:2.22}
  Recall that if a point $c$ is a regular value of a smooth map $f:M\to
  N$ between two manifolds then $(Tf)\inv (c,0) = T(f\inv (c))$.
  Consequently  the tangent  functor $T$ preserves
transverse fiber products.  Namely if  $f:M\to L$, $g:N\to L$ are two maps so
that $(f,g):M\times N\to L\times L$ is transverse to the diagonal
$\Delta_L\subset L\times L$ then
\[
T(M\times_{f,L,g} N)   = T \left( (f,g)\inv (\Delta_L) \right) =
(T(f,g))\inv (T\Delta_L) =  T M\times_{f,TL,g} TN .
  \]
\end{remark}
The following result is well-known.
\begin{lemma}\label{lem:2.23'}
The tangent functor $T:\Man \to \Man$ and the natural transformation
$\pi:T\Rightarrow \id_M$ extend to a functor $T^\LieGpd:\LieGpd\to
\LieGpd$ and to a natural transformation $\pi: T^\LieGpd \Rightarrow \id_\LieGpd$.
 \end{lemma} 

\begin{proof}[Sketch of proof]
To define $T^\LieGpd$ on a Lie groupoid $G$ we apply the the tangent
functor $T:\Man \to \Man$ to all the structure maps of the groupoid
$G$.    Since $T$ is a functor we get a Lie groupoid, the tangent
groupoid of $G$.  This groupoid is commonly
denoted by $TG$.  We set $T^\LieGpd(G):= TG$.  This defines the
functor $T^\LieGpd$ on objects.

A functor $f:G\to H$ between two Lie groupoids consists of a pair of
smooth maps: $f_0:G_0\to H_0$ and $f_1:G_1\to H_1$.  These maps induce a
smooth map $f_2: G_1\times_{s,G_0, t}G_1\to H_1\times_{s,H_0, t} H_1$
and 
intertwine all the structure maps of the two groupoids.  Apply the
functor $T$ gives us smooth maps $Tf_0, Tf_1$ and $Tf_2$ that define a
functor $Tf:TG\to TH$.  This defines $T^\LieGpd$ on 1-cells (i.e., on
1-morphisms).

A natural transformation $\alpha:f\Rightarrow k$ between two functors
$f,k:G\to H$ is a smooth map $\alpha:G_0 \to H_1$ (which is subject to
the appropriated conditions).  Applying the functor $T$ gives us
$T\alpha:TG_0\to TH_1$.  The map $T\alpha$ is a natural transformation from
$Tf$ to $Tk$.  This defines $T^\LieGpd$ on 2-cells.

Given a Lie groupoid $G$ we have two projections $(\pi_G)_0 :TG_0\to
G_0$ and $(\pi_G)_1 :TG_1\to G_1$.  The projections assemble into a
functor $\pi_G:TG\to G$.   It is not hard to check that for any map
$f:G\to H$ of Lie groupoids the diagram
\[
\xy
(-14,10)*+{TG}="1";
(14,10)*+{TH}="2";
(-14, -10)*+{G}="3";
(14, -10)*+{H}="4";
{\ar@{->}^{Tf} "1"; "2"};
{\ar@{->}_{\pi_G} "1"; "3"};
{\ar@{->}^{\pi_H} "2"; "4"};
{\ar@{->}_{f} "3"; "4"};
\endxy
\]
commutes.  Hence the collection of functors $\{\pi_G\}_{G\in \LieGpd}$
defines a natural transformation $\pi:T^\LieGpd \Rightarrow \id$.
\end{proof}

\begin{lemma}\label{lem:2.24'}
The functor $T:\Man\to \Man$ extends to a functor $T^\Bi:\Bi\to \Bi$
from the bicategory $\Bi$ of Lie groupoids, bibundles and isomorphisms
of bibundles to itself.
 \end{lemma} 

\begin{proof}[Sketch of proof]  For a Lie groupoid $G$ we set $T^\Bi
  (G) = TG$.  This defines $T^\Bi$ on 0-cells.   Given a bibundle
  $P:G\to H$ the application of the tangent functor $T:\Man\to \Man$
  gives us the bibundle $TP:TG\to TH$.  Given two bibundles $P,Q:G\to
  H$ and an isomorphism $\alpha:P\rightarrow Q$ of bibundles its
  derivative $T\alpha:TP\to TQ$ is also an isomorphism of bibundles.
  It is not hard to check that for any two Lie groupoids $G$ and $H$
  the map 
  \[
T: \Hom_\Bi (G,H) \to \Hom_\Bi (TG, TH)
    \]
defined above is a functor.  Given a Lie groupoid $G$ we defined the
identity bibundle $\langle \id_G\rangle$ to be the manifold $G_1$
together with the source and target maps as left and right anchors and
left and right multiplications as left and right actions of $G$ on
$G_1$.  Then $T\langle \id_G\rangle = TG_1 = \langle \id_{TG}
\rangle$. Hence the comparison 2-cells $\mu_G :\id_{T^\Bi (G)} \to
T^\Bi (\id_G)$ are identity 2-cells.   We also need the comparison 2-cells
\[
\mu_{Q,P}: T^\Bi Q \circ T^\Bi P \to T^\Bi (Q\circ P).
\]
They are constructed as follows.  Given a pair $G\xrightarrow{P}H
\xrightarrow{Q} K$ of composable bibundles $T(P\times_{H_0} Q) = TQ
\times_{TH_0} TQ$ (cf.\ Remark~\ref{rmrk:2.22}).  Additionally $T (P\times_{H_0} Q)/H)
\simeq\left(T (P\times_{H_0} Q)\right)/TH$ since for any Lie groupoid $H$
 and any $H$-principal bundle $R\to B$, $TR/TH$ is isomorphic to $
 TB$.  Consequently
 \[
T(Q\circ P) = T ((P\times_{H_0} Q)/H) \simeq T(P\times_{H_0} Q)/TH
\simeq (TP\times_{TH_0} TQ)/TH  = TQ\circ TP.
   \]
This diffeomorphism is the desired invertible 2-cell $\mu_{Q,P} $.
This concludes our construction of data for the pseudo-functor $T^\Bi$
(cf.\ Definition~\ref{def:App3}).

It is an unilluminating exercise to check that the data constructed
above satisfies the conditions of Definition~\ref{def:App3}.  We leave
it to an interested reader.
\end{proof}

\begin{lemma}\label{lem:2.25'}
The functors $U\circ T^\LieGpd$ and $T^\Bi \circ
U$ are isomorphic.  Here $U:= \langle \,
\rangle : \LieGpd \to \Bi$ is the localization functor. 
  \end{lemma}

  \begin{proof}
We construct a pseudo-natural transformation $\sigma: U\circ
T^\LieGpd\Rightarrow T^\Bi \circ U$ as follows (cf.\
Definition~\ref{def:App_nat}).  Since for a groupoid $G$
\[
U\circ T^\LieGpd (G) = TG = T^\Bi\circ U (G)
  \]
we set $\sigma_G : =\langle \id_{TG}\rangle$.  Given a functor $f:G\to
H$ we define the 2-cell $\sigma_f: \sigma_H \circ T^\Bi (U(f))
\Rightarrow U (T^\LieGpd (f)) \circ \sigma_G$ as composite of the
isomorphisms of bibundles
\[
\langle \id_{TH} \rangle \circ T\langle f\rangle \to T\langle f\rangle
\to \langle Tf \rangle \to \langle Tf \rangle \circ \langle \id_{TG}\rangle.
\]
  Here the middle arrow is the diffeomorphism $T(G_0\times_{f,H_0,
    t}H_1) \xrightarrow{\simeq} TG_0 \times_{Tf, TH_0, Tt} TH_1$.
Given an isomorphism $\alpha: f\Rightarrow f'$ between two functors
$f,f':G\to H$ the diagram
\[
\xy
(-20,10)*+{T(G_0\times_{f,H_0, t}H_1)}="1";
(24,10)*+{TG_0 \times_{Tf, TH_0, Tt} TH_1}="2";
(-20, -10)*+{T(G_0\times_{f',H_0, t}H_1)}="3";
(24, -10)*+{TG_0 \times_{Tf', TH_0, Tt} TH_1}="4";
{\ar@{->}^{\sim} "1"; "2"};
{\ar@{->}_{T\langle \alpha \rangle} "1"; "3"};
{\ar@{->}^{\langle T\alpha \rangle} "2"; "4"};
{\ar@{->}_{\sim} "3"; "4"};
\endxy
\]
commutes.   It follows that $\sigma_f$'s are components of a natural
transformation
\[
\sigma_{G,H}: \left(U\circ T^\LieGpd\right) \circ \sigma_G \Rightarrow
\sigma_H \circ \left( T^\Bi \circ U\right).
\]
It is easy to see that with these definitions a version of
\eqref{eq:App2} commutes.  Checking that \eqref{eq:App1} commutes  is a bit
harder. It amounts to proving the following fact: given three Lie
groupoids $A$, $B$ and $C$ and two functors $A\xrightarrow{f}B$,
$B\xrightarrow{g}C$ the two isomorphisms of the bibundles
\[
(T\langle g\rangle \circ  T\langle f\rangle)\circ \langle \id_A\rangle
\xrightarrow{ (\tau \star \id_{\sigma_A})} (T \langle g\circ f
\rangle)\circ \langle \id_{TA}\rangle \xrightarrow{\sigma_{gf} }
\langle \id _{TC} \circ \langle T(gf)\rangle
\]
and
\begin{gather*}
\left( T\langle g\rangle \circ T\langle f\rangle\right) \circ \langle
\id_{TA}\rangle \to T\langle g\rangle \circ \left(
  T\langle f\rangle \circ \langle \id_{TA} \right) \xrightarrow{\id
  _{T\langle g\rangle}\star \sigma_f } T\langle g\rangle \circ \left(
  \id_{TB} \circ \langle Tf\rangle\right) \to \\
\left( \langle g\rangle \circ \langle \id_{TB}\rangle\right) \circ
\langle Tf \rangle \to \left(\langle \id_{TC}\rangle \circ \langle
  Tg\rangle \right) \circ \langle Tf \rangle \to \langle
\id_{TC}\rangle \circ \left( \langle Tg\rangle \circ \langle Tf\rangle
\right) \to \langle \id _{TC} \circ \langle T(gf)\rangle
\end{gather*}  
 are equal.  
  \end{proof}

\begin{lemma} \label{lem:2.26'}
There exists a transformation $\tilde{\pi}: T^\Bi \Rightarrow
\id_\Bi$ whose 1-cells $\tilde{\pi}_G:T^\Bi (G) \to G$ are (isomorphic
to) the bibundles $\langle \pi_G\rangle$, where as before $\pi_G:TG\to
G$ are the projection functors.
\end{lemma}

\begin{proof}
Since $U  = \langle \, \rangle :\LieGpd \to \Bi$ is a localization
functor, the pullback by $U$ defines an equivalence of bicategories
\[
U^*:= -\circ U: \Hom(\Bi,\Bi) \to \Hom_W (\LieGpd, \Bi).
\]
Here as before $\Hom_W (\LieGpd, \Bi)$ denotes the bicategory of
functors that send essential equivalence in $\LieGpd$ to invertible
bibundles.  Consequently for any two functors $F,G:\Bi \to \Bi$ the
functor
\[
U^*: \Hom(F,G) \to \Hom (F\circ U, G\circ U)
  \]
is an equivalence of categories. Note that the objects of $\Hom(F,G)$
are pseudo-natural transformations and morphisms are modifications.
In Lemma~\ref{lem:2.23'} we constructed a natural transformation
$\pi:T^\LieGpd \to \id$ and Lemma~\ref{lem:2.25'} we constructed a
natural isomorphism $\sigma:U\circ T^\LieGpd \Rightarrow T^\Bi \circ
U$.  Therefore we have a natural transformation $U\pi \circ \sigma\inv
: T^\Bi \circ U\Rightarrow U$.    Since $U^*: \Hom(T^\Bi, \id)\to
\Hom(T^\Bi \circ U, U)$ is essentially surjective, there exists a
pseudo-natural transformation $\tilde{\pi}: T^\Bi \to \id$ so that
$\tilde{\pi}\circ U$ differs from $U\pi\circ \sigma\inv$ by a
modification.   It will be convenient to fix one such modification
throughout the paper.  It follows that for each Lie groupoid $G$  we
have chosen an isomorphism of bibundles $\tilde{\pi}_G \to \langle
\pi_G\rangle$ where  $\tilde{\pi}_G$ is the component of the
transformation  $\tilde{\pi}$ at $G$.
\end{proof}

\section{The Lie 2-algebra $\X(G)$ of vector fields on a Lie
  groupoid $G$}
\label{sec:3}
In this section we prove Theorem~\ref{lem:03.4}: the category of
multiplicative vector fields on a Lie groupoid underlies a strict
Lie 2-algebra.  We start by recalling the definition of the category
of multiplicative vector fields.

As we saw in Subsection~\ref{subsec:tangent} for any Lie groupoid $G$
we have the tangent groupoid $TG$ and a functor $\pi_G:TG\to G$.

\begin{definition}[Hepworth \cite{Hepworth}]  
\label{def:XG}
  Consider a Lie groupoid $G$ with its tangent groupoid $\pi_G:TG \to
  G$.  The {\sf category $\X(G)$ of multiplicative vector fields} is
  defined as follows.  The set of {\sf objects} of $\X(G)$ is
\[
\X(G)_0:= \{v:G\to TG\mid v\textrm{ is a functor and } \pi_G\circ v = \id_G\}.
\]
This is the set of multiplicative vector fields of Mackenzie and
Xu~\cite{MackXu}.  A {\sf morphism} in $\X(G)$ from a multiplicative vector
field $v$ to a multiplicative vector field $w$ is a natural
transformation $\alpha:v\Rightarrow w$ such that $\pi_G \star \alpha =
1_{id_G}$.  That is, for every point $x\in G_0$ we require that
\begin{equation} \label{eq:xg-mor}
\pi_G (\alpha_x) = 1_x.
\end{equation}
The composition of morphisms is the vertical composition of natural
transformations.  Note that every morphism of $\X(G)$ is automatically
invertible (since $TG$ is a groupoid). Hence the category $\X(G)$ is a groupoid.
\end{definition}

\begin{notation}\label{not:3.2}
  We denote the source and target maps in the category $\X(G)$ by
  $\fs$ and $\ft$, respectively.  The unit map is denoted by
  $\mathbf{1}$, the inversion by $(\,)\inv$ and the
  composition/multiplication of morphisms by $\circ$.
\end{notation}

\begin{lemma}\label{lem:2vect}
  The category of multiplicative vector fields $\X(G)$ on a Lie
  groupoid $G$ is a 2-vector space.
\end{lemma}

\begin{proof}
  Mackenzie and Xu proved that the set $\X(G)_0$ of multiplicative
  vector fields is a real vector space  \cite{MackXu}. 

  We next argue that the set of morphisms $\X(G)_1$ of the
  category $\X(G)$ is a vector space as well.  Suppose
  $\alpha_1:v_1\Rightarrow w_1$ and $\alpha_2:v_2\Rightarrow w_2$ are 
  morphisms between multiplicative vector fields.  Equation
  \eqref{eq:xg-mor} says that $\alpha_1$ and $\alpha_2$ are both
  sections of the vector bundle
\[
TG_1|_{G_0}\to G_0
\]
where we have suppressed the unit map $1_G:G_0 \to G_1$.  Clearly the
linear combination $\lambda_1\alpha_1 + \lambda_2\alpha_2$ is again a
section of the bundle $TG_1|_{G_0}\to G_0$ for any choice of scalars
$\lambda_1, \lambda_2$.  We need to check that it is actually a
natural transformation from $\lambda_1v_1 + \lambda_2v_2$ to
$\lambda_1w_1 + \lambda_2w_2$.  That is, we need to check that for any
arrow $y\xleftarrow{\gamma}x$ in the groupoid $G$
\[
(\lambda_1\alpha_1 + \lambda_2\alpha_2)_y \bullet 
((\lambda_1v_1 + \lambda_2v_2)(\gamma)) = 
((\lambda_1w_1 + \lambda_2w_2)(\gamma))\bullet 
(\lambda_1\alpha_1 + \lambda_2\alpha_2)_x.
\]
Here and below $\bullet:TG_1\times_{TG_0}TG_1\to TG_1$ denotes the
multiplication in the Lie groupoid $TG$.

Since $\bullet$ is the derivative of the multiplication
$m:G_1\times_{G_0}G_1\to G_1$ in the groupoid $G$, it is fiberwise
linear: for any $(\gamma_2, \gamma_1)\in G_1\times_{G_0}G_1$ and $(a_1,a_2), (b_1, b_2) \in T_{\gamma_2}
G_1\times_{TG_0}T_{\gamma_1}G_1 = T_{(\gamma_1, \gamma_2)}(
G_1\times_{G_0}G_1)$ we have (in the prefix notation)
\begin{equation}\label{eq:prefix}
\bullet (\lambda (a_1, a_2) + \mu (b_1, b_2) ) = \lambda (\bullet(a_2, a_1)) +
\mu (\bullet (b_1,b_2))
\end{equation}
for all scalars $\lambda, \mu$.
In the infix notation \eqref{eq:prefix} reads
\begin{equation}\label{eq:3.3}
(\lambda a_1 + \mu b_1) \bullet (\lambda a_2 + \mu b_2)=
\lambda (a_1\bullet a_2) + \mu (b_1\bullet b_2).
\end{equation}
Hence 
\begin{eqnarray*}
(\lambda_1\alpha_1 + \lambda_2\alpha_2)_y \bullet 
((\lambda_1v_1 + \lambda_2v_2)(\gamma)) &= & 
\lambda_1 ((\alpha_1)_y \bullet v_1 (\gamma)) 
+ \lambda_2 ((\alpha_2)_y \bullet v_2 (\gamma))\\
&=& \lambda_1 (w_1(\gamma)\bullet (\alpha_1)_x) 
+ \lambda_2 (w_2(\gamma)\bullet (\alpha_2)_x)\\
&=& ((\lambda_1w_1 + \lambda_2w_2)(\gamma))\bullet 
(\lambda_1\alpha_1 + \lambda_2\alpha_2)_x.
\end{eqnarray*}
Here the first and third equalities hold by \eqref{eq:3.3}.  In the
second equality we used the fact that $\alpha_1:v_1\Rightarrow w_1$,
and $\alpha_2:v_2\Rightarrow w_2$ are natural transformations.
Therefore the space of morphisms $\X(G)_1$ is a vector space.

Moreover the computation above shows that for $\lambda_1, \lambda_2\in
\R$, $\alpha_1:v_1\Rightarrow w_1$, $\alpha_2:v_2\Rightarrow w_2\in
\X(G)_1$ the source of $\lambda_1 \alpha_1 + \lambda_2 \alpha_2$ is
$\lambda_1 v_1 + \lambda_2 v_2$.  That is, the source map $\fs:
\X(G)_1\to \X(G)_0$ of the category $\X(G)$ is linear.
Similarly the target map $\ft$ is linear.  It is also
easy to see that the unit map $\X(G)_0\to \X(G)_1$ is linear as well.

Finally we need to check that multiplication/composition $\circ$ in the
category $\X(G)$, which is the vertical composition of natural
transformations, is linear as a map from $\X(G)_1 \times _{\X(G)_0}
\X(G)_1 $ to $\X(G)_1$.  That is, we need to check that 
\begin{equation}\label{eq:third}
(\lambda \alpha_2 + \mu \beta_2)\circ (\lambda \alpha_1 + \mu \beta_1)=
\lambda (\alpha_2\circ \alpha_1) + \mu (\beta_2 \circ \beta_1)
\end{equation}
for all $\lambda, \mu \in \R$, $(\alpha_2, \alpha_1), (\beta_2,
\beta_1)\in \X(G)_1 \times _{\X(G)_0} \X(G)_1$.  Recall that the
vertical composition $\circ$ is computed pointwise: for any
composable natural transformations $\delta_1, \delta_1$ and any point
$x\in G_0$
\[
(\delta_2 \circ \delta_1)_x = (\delta_2)_x \bullet (\delta_1)_x,
\]
where as before $\bullet$ is the multiplication in $TG$.  Since
$\bullet$ is fiberwise linear \eqref{eq:third} follows.  This
concludes our proof that the category $\X(G)$ of vector fields on a
Lie groupoid $G$ is internal to the category of vector spaces, that
is, $\X(G)$ is a 2-vector space.
\end{proof}

\begin{theorem} \label{lem:03.4}
The category of vector fields $\X(G)$ on a Lie groupoid
$G$ is a (strict) Lie 2-algebra.
\end{theorem}

\begin{proof}
  Recall the notation: $s:G_1\to G_0$ is the source map for the
  groupoid $G$, its differential $Ts:TG_1\to TG_0$ is the source map
  for the tangent groupoid $TG$.  We use $\fs, \ft$ to denote the
  source and target maps of the groupoid $\X(G)$, respectively.

  By Lemma~\ref{lem:Lie2_from_cross} it is enough to: (i) give the vector
  spaces $\ker (\fs: \X(G)_1\to \X(G)_0)$ and $\X(G)_0$ the structure
  of Lie algebras, (ii) check that $\partial:= \ft|_{\ker \fs}:\ker\fs \to 
  \X(G)_0$ is a Lie algebra map, (iii) define an action $D:\X(G)_0\to \Der
  (\ker \fs))$ on $\ker \fs$ by derivations, and (iv) check the
  compatibility of $\partial$ and $D$:
\begin{eqnarray}
\partial(D(X) \alpha) &=& [X, \partial(\alpha)], \label{eq:5.1nv}\\
D(\partial \alpha_1) \alpha_2 &=& [\alpha_1, \alpha_2]\label{eq:5.2nv}
\end{eqnarray}
for all $\alpha, \alpha_1, \alpha_2 \in \ker \fs$ and all
multiplicative vector fields $X$ on the Lie groupoid $G$
(compare with  Definition~\ref{def:cross}).

The fact that the vector space $\X(G)_0$ of multiplicative vector
fields carries a Lie bracket is due to Mackenzie and Xu \cite{MackXu}.
  We  argue next that $\ker \fs$ is the space of sections of the Lie
  algebroid $A_G\to G_0$.  By definition of the source map $\fs$,
\[
\ker \fs =\{ \alpha:X\Rightarrow Y \mid X=0\}.
\]
Therefore $\alpha \in \ker \fs$ if and only if there is a
multiplicative vector field $Y$ so that the diagram 
\begin{equation}\label{eq:3.5nv}
\xy
(-10,10)*+{0_x}="1";
(20,10)*+{0_y}="2";
(-10, -10)*+{Y(x)}="3";
(20, -10)*+{Y(y)}="4";
{\ar@{->}^{0_\gamma} "1"; "2"};
{\ar@{->}_{\alpha_x} "1"; "3"};
{\ar@{->}^{\alpha_y} "2"; "4"};
{\ar@{->}_{Y(\gamma)} "3"; "4"};
\endxy
\end{equation}
commutes for all arrows $y\xleftarrow{\gamma}x$ in $G_1$.  Hence if
$\alpha\in \ker \fs$ then $Ts(\alpha_x) = 0_x$ for all $x\in G_0$.
That is, $\alpha$ is a section of $A_G\to G_0$.  Conversely if
$\alpha:G_0\to A_G$ is a section of the Lie algebroid we can define a
multiplicative vector field $Y:G\to TG$ so that \eqref{eq:3.5nv}
commutes.  Namely on objects we define 
\[
Y(x):= Tt (\alpha_x) \qquad \textrm{for all  } x\in G_0.
\]
And for  $y\xleftarrow{\gamma}x$ in $G_1$ we set
\begin{equation}\label{eq:3.6vn}
Y(\gamma) = \alpha_y \bullet 0_\gamma\bullet (\alpha_x)\inv.
\end{equation}
Here as before $\bullet$ is the multiplication in $TG$ and $(
\,)\inv =Ti:TG_1\to TG_1$ is the inverse map, which is the derivative of
the inverse map $i$ of the groupoid $G$.  We conclude that 
\[
\ker ( \fs: \X(G)_1\to \X(G)_0) = \Gamma (A_G) 
\]
and that
\[
\partial := \ft |_{\ker \fs}:\ker \fs \to \X(G)_0 
\]
is given by
\begin{equation}\label{eq:3.7nv}
(\partial \alpha)(\gamma) = 
\alpha_{t(\gamma)}\bullet 0_\gamma\bullet (\alpha_{s(\gamma)})\inv
\end{equation}
for all $\gamma \in G_1$. Note that \eqref{eq:3.7nv} can be written as 
\begin{equation}\label{eq:3.8nv}
\partial \alpha= \overrightarrow{\alpha}+\overleftarrow{\alpha}.
\end{equation}
where 
\[
\overrightarrow{\alpha}(\gamma) = TR_\gamma\, \alpha (t(\gamma))
\]
and 
\[
\overleftarrow{\alpha}(\gamma) = T(L_\gamma\circ i)\, \alpha (t(\gamma))
\]
for all $\gamma\in G_1$.  Here $R_\gamma$ and $L_\gamma$ are right and
left multiplications by $\gamma$, respectively. 
The bracket on the space of sections $\Gamma(A_G)$ of the
Lie algebroid $A_G\to G_0$ is defined by requiring that the injective
map
\[
^\rightarrow:\Gamma(A_G)\to \Gamma(TG_1), 
\qquad \alpha \mapsto \overrightarrow{\alpha}.
\]
is a map of Lie algebras.   Consequently 
\[
^\leftarrow:\Gamma(A_G)\to \Gamma(TG_1), \qquad \alpha \mapsto
\overleftarrow{\alpha}
\]
is also a map of Lie algebras.   Since left -and right-invariant vector fields commute (cf.\ \cite{MacK} ) we conclude that
\[
\partial = \ft|_{\ker \fs}:\ker\fs = \Gamma(A_G) \to 
  \X(G)_0
\]
is a Lie algebra map.

Following Mackenzie and Xu we define the map
$D$ from the space $\X(G)_0$ of multiplicative vector fields to
$\Hom(\Gamma (TG_1|_{G_0}), \Gamma (TG_1|_{G_0}))$ by setting
\[
D(X)\alpha:= [X, \overrightarrow{\alpha}]|_{G_0}
\]
for all multiplicative vector fields $X$ and all sections $\alpha \in
\Gamma (A_G)$.  Mackenzie and Xu prove \cite[Proposition~3.7]{MackXu}
that $[X, \overrightarrow{\alpha}]$ is tangent to the fibers of $s$
and is right invariant. Hence $[X, \overrightarrow{\alpha}]|_{G_0}$ is
a section of the Lie algebroid $A_G\to G_0$.  They furthermore show
\cite[Proposition~3.8]{MackXu} that $D(X)$ is a derivation of
$\Gamma(A_G)$ and that $D:\X(G)_0 \to \Der(\Gamma(A_G)) $ is a map of
Lie algebras.

Since left- and right-invariant vector fields commute, for any
$\alpha_1, \alpha_2\in \Gamma(A_G)$ we have
\[ [\partial \alpha_1, \overrightarrow{\alpha}_2] =[
\overrightarrow{\alpha}_1+ \overleftarrow{\alpha_1},
\overrightarrow{\alpha}_2] = [\overrightarrow{\alpha}_1,
\overrightarrow{\alpha}_2]
\]
and \eqref{eq:5.2nv} follows.  

We end the proof by showing that \eqref{eq:5.1nv} holds.  On the right
hand side we have
\[
[X, \partial \alpha] = [X, \overrightarrow{\alpha} +  
\overleftarrow{\alpha}]=[X,\overrightarrow{\alpha}]+[X,\overleftarrow{\alpha}]
\]
while on the left,
\[
\partial \left(D(X)\, \alpha\right) =  \left(D(X)\, \alpha\right)^{\to}
+ \left(D(X)\, \alpha\right)^{\leftarrow}.
\]
By definition of $D$,  
\[
\left(D(X)\, \alpha\right)^{\to}= [X,\overrightarrow{\alpha}],
\]
so it remains to prove that $[X, \overleftarrow{\alpha}] = \left(D(X)\,
  \alpha\right)^{\leftarrow}$.  Since $X$ is a functor, 
\[
Ti \circ X = X \circ i.
\]
The inversion map $i$ relates  right- and left-invariant vector fields. That is,
\[
Ti \circ  \overrightarrow{\alpha} = \overleftarrow{\alpha} \circ i
\]
for all $\alpha$.   Consequently 
\[
\left(D(X)\, \alpha\right)^{\leftarrow} (g) = T(L_g \circ i)(D(X)\,
\alpha)(1_{s(g)})= T(L_g) Ti ([X, \overrightarrow{\alpha}] (1_{s(g)})= 
TL_g [X, \overleftarrow{\alpha}](i(1_{s(g)})).
\]
Since $[X, \overleftarrow{\alpha}]$ is left-invariant, 
$TL_g [X, \overleftarrow{\alpha}](i(1_{s(g)})) 
= [X, \overleftarrow{\alpha}] (g)$.  Therefore,
\[
\left(D(X)\, \alpha\right)^{\leftarrow} (g)= [X, \overleftarrow{\alpha}] (g)
\]
for all $g\in G_1$ and we are done.

\end{proof}

\section{Morita invariance of the Lie 2-algebra of vector
  fields %
}
 \label{sec:4}

The goal of this section is to prove 

\begin{theorem} \label{thm:04.1}
The assignment $G\mapsto \X(G)$ of the category of vector fields to a Lie groupoid 
extends to a functor
\begin{equation}\label{eq:Xfunct}
\X: \Bi_\iso \to \Lietwo
\end{equation}
from the bicategory $\Bi_\iso$ of Lie groupoids, invertible bibundles
and isomorphisms of bibundles to the bicategory $\Lietwo$ of Lie
2-algebras.  Hence, in particular, if $P:G\to H$ is a
Morita equivalence of Lie groupoids then $\X(P):
\X(G)\to \X(H)$ is a (weakly) invertible 1-morphism of
Lie 2-algebras in the bicategory $\Lietwo$.
\end{theorem}

Our strategy for constructing the functor $\X$ is to first construct it on
a simpler category.

\begin{definition}
  An {\sf essentially surjective open embedding} %
  of Lie groupoids is a
  functor $f:U\to G$ so that 
\begin{enumerate}
\item The maps on objects $f_0:U_0 \to G_0$ and on morphisms
  $f_1:U_1\to G_1$ are open embeddings and
\item the functor $f$ is an essential equivalence, i.e., the corresponding bibundle 
\[
\langle f \rangle := U_0\times _{f_0, G_0, t}G_1: U\to G
\]
is weakly invertible.  Equivalently $\langle a^R: \langle f \rangle
\to G_0$ is a principal $U$-bundle.
\end{enumerate}
\end{definition}

It is clear that the identity functors are essentially surjective open
embeddings.  Moreover the composition of essentially surjective
open embeddings is again an essentially surjective open embedding.
Consequently Lie groupoids, essentially surjective open embeddings and
natural transformations form a 2-category.

\begin{notation}
  We denote the 2-category of Lie groupoids, essentially surjective open
  embeddings and natural transformations by $\ESOE$.
\end{notation}

\begin{theorem}\label{thm:loc}
  The localization of the bicategory $\ESOE$ at the class $W$ of all
  1-morphisms is the bicategory $\Bi_\iso$ of bicategory of Lie
  groupoids, invertible  bibundles and isomorphisms of
  bibundles.
\end{theorem}

\begin{proof} We apply Proposition~\ref{prop:Pronk}.  Consider the
  localization functor $\langle \,\,\rangle : \LieGpd \to \Bi$
  introduced in Remark~\ref{rmrk:02.2}.  By definition of the
  2-category $\ESOE$ the restriction of the functor $\langle
  \,\,\rangle$ to $\ESOE$ sends every 1-morphism $w:U\to G$ of $\ESOE$
  to an invertible bibundle $\langle w \rangle$ (and a 2-morphism to
  an isomorphism of bibundles).  This gives us a functor
\begin{equation}\label{eq:loc}
\langle \,\, \rangle :\ESOE \to \Bi_\iso, \qquad 
(G \xrightarrow{w}H)\mapsto (G\stackrel{\langle w \rangle}{\to} H).
\end{equation}
The functor is surjective on objects.  By Theorem~\ref{thm:folk1} the
functor is fully faithful on 2-morphisms.

It remains to check that given an invertible bibundle $P:G\to H$ there
exist essentially surjective open embeddings $w_G$, $w_H$ so that
$\langle w_H\rangle \circ P$ is isomorphic to $\langle w_G \rangle$.
Since the bibundle $P$ is weakly invertible, it gives rise to the
\emph{linking groupoid} \cite[Proposition~4.3]{Wei}, denoted $G*_PH$
and recalled presently. The manifold of objects $(G*_PH)_0$ is the
disjoint union $G_0\sqcup H_0$ of the objects of the groupoids $G$ and
$H$.  The manifold of arrows $(G*_PH)_1$ is the disjoint union
$G_1\sqcup P\sqcup P^{-1}\sqcup H_1$.  We think of the manifold $P$ as
the space of arrows from the points of $H_0$ to the points of $G_0$.
We think of the elements of $P\inv$ as the inverses of the elements of
$P$.  The multiplication in $G*_PH$ comes from the multiplications in
the groupoids $G$ and $H$ and the actions of $G$ and $H$ on $P$ and on
$P\inv$.  The inclusion $w_G:G\to G*_PH$ is given by the open
embeddings
\[
G_0\hookrightarrow G_0\sqcup H_0, \qquad G_1 \hookrightarrow G_1\sqcup
P\sqcup P^{-1}\sqcup H_1.
\] 
It is easy to see that $w_G$ is an essential equivalence, i.e., that
the bibundle $\langle w_G\rangle$ is biprincipal, hence weakly
invertible.  Similarly we have the essentially surjective open
embedding
\[
w_H:H\hookrightarrow G*_PH.
\]
A computation shows that the bibundles $\langle w_H\rangle \circ P$ 
and $\langle w_G \rangle$ are isomorphic.
\end{proof}

\begin{proposition} \label{prop:X_on_ESOP}
The assignment
\[
G\mapsto \X(G)
\]
of the Lie 2-algebra of vector fields to a Lie groupoid extends to a
contravariant functor
\begin{equation}\label{eq:w*}
 (\ESOE)^{op} \to \Lietwo_{strict}, \quad (G\xrightarrow{w}  H) 
\mapsto (\X(H)\xrightarrow{w^*}  \X(G))
\end{equation}
from the bicategory $\ESOE$ of Lie groupoids, essentially surjective
open embeddings and natural isomorphism to the strict 2-category
$\Lietwo_{strict}$ of Lie 2-algebras.
\end{proposition}
\begin{proof}
Consider an essentially surjective
open embedding $w:G\to H$.
Then $w(G)\subset H$ is an open Lie subgroupoid and $w:G\to w(G)$ is
an isomorphism of Lie groupoids.  We now assume without any loss of
generality that $G$ is an open subgroupoid of $H$.  Then the tangent
bundle $TG$ is an open subgroupoid of $TH$.  Moreover, any
multiplicative vector field $v:H\to TH$ restricts to a multiplicative
vector field $v|_G :G\to TG$.  Similarly, a morphism
$\alpha:v\Rightarrow u$ of multiplicative vector fields restricts to a
morphism $\alpha|_G:v|_G \Rightarrow u|_G$.  This gives us a functor
\begin{equation} \label{eq:04.2}
w^*: \X(H)\to \X(G), \qquad 
w^*( \alpha:v\Rightarrow u) = (\alpha|_G:v|_G \Rightarrow u|_G).
\end{equation}
The restriction to an open subgroupoid is a map of 2-vector spaces and
preserves the brackets.  Hence \eqref{eq:04.2} is a map of Lie
2-algebras.
\end{proof}

\begin{definition}[The category $\X_\gen(G)$ of {\sf generalized vector fields
  on a Lie groupoid }$G$]\label{ex:Xgen} \mbox{}\\
  Recall that there is a natural transformation $\tilde{\pi}:T^\Bi
  \Rightarrow \id _\Bi$ (Lemmas~\ref{lem:2.24'} and \ref{lem:2.26'}).
  
An object of the category of generalized vector fields $\X_\gen(G)$
on a Lie
  groupoid $G$ is a pair $(P,\alpha_P)$ where $P:G\to TG$ is a
  bibundle and $\alpha_P: \tilde{\pi}_G \circ P \Rightarrow
  \langle {id_G}\rangle)$ an isomorphism of bibundles.

  A morphism $\beta$ in $\X_\gen(G)$ from
  $(P,\alpha_P)$ to $(Q,\alpha_Q)$ is a map of bibundles $\beta:P\Rightarrow
  Q$ so that
\[
\alpha_{Q} =  \alpha_P \circ (\tilde{\pi}_G \star \beta).
\]
Here as before $\star$ denotes whiskering in $\Bi$, and $\circ$ is the
composition of isomorphisms of bibundles.
\end{definition}

\begin{lemma} \label{lemma:equiv_of_vf2} A weakly
  invertible bibundle $P:G\to H$ between two Lie groupoids 
  induces an equivalence of categories
\[
 P_*: \X_\gen(G)\to \X_\gen(H)
\]
between the corresponding categories of generalized vector fields.
\end{lemma}

\begin{proof}
 Since the 1-morphism $P$ is (weakly) invertible, there is 2-morphism
\[
\gamma: (P\circ \langle \id_G \rangle)\circ P\inv \Rightarrow \langle \id_H\rangle .
\]
Given an object $(X,\alpha_X)$ of $\X_\gen (G)$ we define
\[
P_*X:=TP\circ (X\circ P^{-1})
\]
The 2-morphism $\alpha_{P_*X}: \tilde{\pi}_H \circ
P_*X\Rightarrow  \langle \id_H\rangle$ comes from the 2-commutative
diagram
\[
\xy
(-50,0)*+{H}="1";
(-30,0)*+{G}="2";
(-10, 0)*+{TG}="3";
(20,0)*+{TH}="4";
(-10,-20)*+{G}="5";
(20,-20)*+{H.}="6";
{\ar@{->}^{P^{-1}} "1"; "2"};
{\ar@{->}^{X} "2"; "3"};
{\ar@{->}^{TP} "3"; "4"};
{\ar@{->}^{\tilde{ \pi}_G} "3"; "5"};
{\ar@{->}^{\tilde{\pi}_H } "4"; "6"};
{\ar@{->}^{P} "5"; "6"};
{\ar@{->}^{P} "5"; "6"};
{\ar@{->}_{\langle \id_G\rangle } "2"; "5"};
{\ar@{<=}^{\tilde{\pi}_P} (2,-12);(6,-8)} ; 
{\ar@{<=}^{\alpha_X} (-17,-8);(-13,-4)} ; 
{\ar@{<=}_{\gamma} (-25,-22);(-21,-18)} ; 
{\ar@/_04.pc/@{->}_{\langle id_H\rangle} "1";"6"} ;
\endxy 
\]
(The 2-morphism $\tilde{\pi}_P$ is part of the data of the natural
transformation $\tilde{\pi}: T^\Bi \Rightarrow \id_\Bi$; see Appendix~\ref{sec:app}.)
We set
\[
\alpha_{P_*X}:= \gamma \circ (P\star \alpha_X \star P\inv )\circ (\tilde{\pi}_P
\star (X\circ P\inv)).
\]
Given a morphism $\beta: (X, \alpha_X) \to (Y, \alpha_Y)$ in the
category $\X_\gen  (G)$ we define
\[
P_*\beta:=TP\star \beta \star P^{-1}.%
\]
A diagram chase ensures that $P_*\beta$ is indeed  a morphism in
$\X_\gen (H)$ from 
$(P_*X, \alpha_{P_*X})$ to $(P_*Y,\alpha_{P_*Y})$.

Finally one checks that the functor $(P\inv)_*: \X(H) \to
\X(G)$ is a weak inverse of $P_*$.  Hence $P_*$ is an equivalence
of categories as claimed.
\end{proof}

\begin{definition} \label{def:incl_functor} The ``inclusion'' functor
  $\imath_G: \X(G)\hookrightarrow \X_\gen(G)$ is defined as follows.
  Suppose $v:G\to TG$ is a multiplicative vector field. Since
  $U: \LieGpd\to \Bi$ is a functor the composition
  $\langle \pi_G \rangle \circ \langle v \rangle $ is isomorphic to
  $\langle \pi_g \circ v\langle = \langle \id_G \rangle$.  Since we
  fixed the modification from $\tilde{\pi}$ to $\pi \circ U$ there is
  a canonical isomorphism of bibundle
  $\tilde{\pi}_G \to \langle \pi_G \rangle$.  Consequently there is a
  canonical isomorphism
  $\alpha_{\langle v \rangle}:\tilde{\pi}_G \circ \langle v \rangle
  \Rightarrow \langle \id_G\rangle$.  We define
  \[
\imath_G (v) := (\langle v \rangle, \alpha _{\langle v \rangle}).
\]
Given a morphism $v\xrightarrow{\gamma} v'$ in $\X(G)$ we define
\[
  \imath_G ({\gamma} ) = \langle \gamma \rangle;
\]
it is a morphism in $\X _\gen(G) $ from $(\langle v \rangle,
\alpha_{\langle v \rangle})$ to $(\langle v' \rangle,
\alpha_{\langle v' \rangle})$. 
  \end{definition}

\begin{theorem}\label{thm:04.3}
For any Lie groupoid $G$ the  inclusion functor 
\[
\imath_G: \X(G)\hookrightarrow \X_\gen(G)
\]
defined above is an equivalence of categories.
\end{theorem}
\begin{remark}
  In the case where the groupoid $G$ is proper Theorem~\ref{thm:04.3}
  can be deduced from \cite[Theorem 15]{Hepworth}; see
 Remark~\ref{rmrk:0proper} below.
\end{remark}

The proof of Theorem~\ref{thm:04.3} is technical; we carry it out in
section~\ref{sec:6} below.

\begin{lemma} \label{lemma:key}
  Let $w:G\to G'$ be an essentially surjective open embedding.  Then
  (the functor of categories underlying) the pull-back/restriction functor
  $w^*:\X(G') \to \X(G)$ 
  which is given by
  \[
    w^* (v\xrightarrow{\gamma} v') = (v|_G \xrightarrow{\gamma|_G} v
'|_G)  
    \]
is fully faithful and essentially
  surjective.
\end{lemma}

\begin{proof}  We first argue that $w^*$ is fully faithful.
We want to show that given a morphism $\delta:v|_G\to v'|_G$ in
$\X(G)$ there exists a unique arrow $\tilde{\delta}:v\to v'$ in
$\X(G')$ so that $\tilde{\delta}|_G = \delta$.

We deal with uniqueness
first.  Suppose $\gamma, \gamma': v\to v'$ are two morphisms in
$\X(G')$ with $\gamma|_G =  \gamma'|G$.   Fix an object $y$ of $G'$. 
Since $w:G\to G'$ is
essentially surjective, for any object $y$ of $G'$ there is an arrow
$\mu:y\to x$ with $x$ an object of $G$.    Then $\gamma_x = \gamma'_x$
and therefore 
\item 
\[
\gamma_y = v'(\mu\inv) \circ \gamma _x \circ v(\mu) = v'(\mu\inv)
\circ \gamma' _x \circ v(\mu) = \gamma'_y.
  \]
  We conclude that $\gamma = \gamma'$.  Now given $\delta:v|_G\to v'|_G$
  we define $\tilde{\delta}:v\to v'$ at an object $y$ by choosing as
  above $\mu:y\to x$ and setting
  \[
\tilde{\delta}_y = v'(\mu\inv) \circ \delta _x \circ v(\mu).
    \]
    If $\nu:y\to x$ is another arrow then
    \[
v'(\nu \circ \mu\inv ) \circ \delta _x  = \delta_x \circ v (\nu \circ \mu\inv ).
\]
Therefore
\[
v'(\mu\inv ) \circ \delta_x \circ v(\mu) = v'(\nu\inv ) \circ \delta_x
\circ v(\nu).
\]
Moreover the dependence of $\tilde{\delta} $ on $y$ is smooth: since
$w:G\to G'$ is an essential equivalence the right anchor $a^R_{\langle
  w\rangle} :\langle w \rangle \to G_0'$ is a surjective submersion.
Note that $\langle w \rangle = t\inv (G_0) \subset G_1'$ and $a^R_{\langle
  w\rangle} (\mu) = s(\mu)$.   Choose a local section $\sigma: U\to
\langle w \rangle$ of $a^R_{\langle w\rangle}$ with $y\in U$.    Then
for all points $z\in U$
\[
\tilde{\delta} (z) = v' ((\sigma (z))\inv ) \circ \delta
_{t(\sigma(y))} \circ v (\sigma (z)),
  \]
  which is smooth.   We conclude that $w^*$ is fully faithful.

  To prove essential surjectivity we need to argue that for any multiplicative
  vector field $u:G\to TG$ there is a multiplicative vector field
  $\tilde{u}: G'\to TG'$ and an isomorphism $\tilde{u}|_G
  \xrightarrow{\sim} u$.  The functor $\imath_{G'} : \X(G')\to
  \X_\gen(G')$ is an equivalence of categories by Theorem~\ref{thm:04.3}.
The functor  $(\langle w \rangle )_*: \X_\gen (G) \to \X_\gen(G')$ is
an equivalence of categories by Lemma~\ref{lemma:equiv_of_vf2}.
Therefore for any $u\in \X(G)_0$ there is a vector field $\tilde{u}
\in \X(G')$ and an isomorphism
\[
\langle \tilde{u}\rangle \xrightarrow{\sim} \langle w\rangle_* \langle
u \rangle.
  \]
  Since
  \[
\xy
(-14,10)*+{TG}="1";
(14,10)*+{TG'}="2";
(-14, -10)*+{G}="3";
(14, -10)*+{G'}="4";
{\ar@{->}^{T w} "1"; "2"};
{\ar@{->}^{\tilde{u}|_G} "3"; "1"};
{\ar@{->}_{\tilde{u}} "4"; "2"};
{\ar@{->}_{w} "3"; "4"};
\endxy
\]
commutes, since $U:\LieGpd \to \Bi$ is a functor and since $T\langle
w\rangle$ is isomorphic to $\langle Tw \rangle$
\[
T\langle w\rangle \circ \langle \tilde{u}|_G \rangle
\xrightarrow{\sim} \rangle \tilde{u}\rangle \circ \langle w\rangle.
  \]
  It follows that $\langle \tilde{u}\rangle $ is isomorphic to
  $T\langle w \rangle \circ (\langle \tilde{u}|_G \rangle \circ
  \langle w \rangle \inv )$.  Therefore
  \[
    T\langle w \rangle \circ (\langle \tilde{u}|_G \rangle \circ
    \langle w \rangle \inv ) \xrightarrow{\sim} T\langle w \rangle
    \circ (\langle u\rangle \circ \langle w \rangle \inv )
    \]
 Consequently there is an isomorphism $\beta: \langle \tilde{u}|_G
 \rangle \to \langle u\rangle $ of bibundles.  Note that this {\em
   does not } yet imply that the generalized  vector fields $(\langle
 \tilde{u}|_G \rangle, \alpha _{\langle
 \tilde{u}|_G \rangle} ) $ and $(\langle
 {u} \rangle, \alpha _{\langle
 u\rangle} ) $ are isomorphic in the category $\X_\gen(G)$.  The issue
is $(\langle
 \tilde{u}|_G \rangle, \alpha _{\langle
 \tilde{u}|_G \rangle} ) $   need not equal $\alpha_{\langle u\rangle}
\circ (\tilde{\pi}_G \star \beta)$.     But
\[
\beta: ( \langle \tilde{u}|_G
 \rangle , \tau ) \to (\langle u\rangle, \alpha _{\langle u\rangle})
  \]
is a morphism in $\X_\gen(G)$ if we set $\tau = \alpha_{\langle u\rangle}
\circ (\tilde{\pi}_G \star \beta)$.  By Lemma~\ref{lem:7.5}  $( \langle \tilde{u}|_G
 \rangle , \tau )$ is isomorphic to $( \langle \tilde{u}|_G
 \rangle , \alpha_{\rangle \tilde{u}|_G \rangle } )$.  Consequently $
\imath_G (\tilde{u}|_G) \textrm{ is isomorphic to } \imath_G (u)$.
Since $\imath_G$ is an equivalence of categories $\tilde{u}|_G$ is
isomorphic to $u$ and we are done.

  \end{proof}

\begin{lemma}\label{prop:X_on_ESOP2}
  The functor \eqref{eq:w*} of Proposition ~\ref{prop:X_on_ESOP} takes
  every essentially surjective open embedding to a weakly invertible
  1-morphism of Lie 2-algebras.
\end{lemma}

\begin{proof}
  By Lemma~\ref{lemma:key} the functor $w^*:\X(G')\to \X(G)$
  associated to an essentially surjective open embedding
  $w:G\to G'$  is fully faithful and essentially surjective, hence an
  essential equivalence of Lie 2-algebras.
  The localization functor $\langle \, \rangle :\Lietwo_{strict}
  \to \Lietwo$ takes all essential equivalences to weakly invertible 1-
  morphisms.
\end{proof}

We are now in position to extend the assignment $G\mapsto
\X(G)$ to a (covariant) functor $\X:\ESOE \to \Lietwo$.
\begin{definition}
We define the functor $\X:\ESOE \to \Lietwo$ on objects to be the assignment 
\[
G\mapsto \X(G).
\]
Given an essentially surjective open embedding $G\xrightarrow{w}G'$,
the bibundle $\langle w^*\rangle$ is weakly invertible in $\Lietwo$ by
Lemma~\ref{prop:X_on_ESOP2}.  We set
\[
\X (w):= (\langle w^* \rangle) \inv.
\]
\end{definition}

\begin{proof}[Proof of Theorem~\ref{thm:04.1}]

  Since $\langle \,\,\rangle :\ESOE \to \Bi_\iso$ is a localization of
  the bicategory $\ESOE$ at the class $\ESOE_1$ of all 1-morphisms and
  since the functor $\X:\ESOE\to \Lietwo$ sends every 1-morphism of
  $\ESOE$ to an invertible morphism there exists by
  Proposition~\ref{prop:Pronk}  functor
\[
\tilde{\X}:\Bi_\iso \to \Lietwo
\]
(which is unique up to isomorphism)
and an isomorphism $\tilde{\X}\circ \langle \,\,\rangle
\stackrel{\sim}{\Leftrightarrow} \X$ of functors.  It is no loss of
generality to assume that $\tilde{\X}(G) = \X(G)$ for every Lie
groupoid $G$.  We now drop the \,$\tilde{\,}$\, and obtain the desired
functor $\X:\Bi_\iso \to \Lietwo$.
\end{proof}

We end the section with a result that we will need in
Section~\ref{sec:5}.

\begin{lemma} \label{lemma:oldkey}
  Let $w:G\to G'$ be an essentially surjective open embedding,
  $w^*:\X(G') \to \X(G)$ the pull-back/restriction functor of
  Proposition~\ref{prop:X_on_ESOP} and $(\langle w \rangle)_*:
  \X_\gen(G)\to \X_\gen (G')$ the push-forward along the bibundle
  $\langle w \rangle$ constructed in Lemma~\ref{lemma:equiv_of_vf2}.
  Then the diagram
\begin{equation} \label{diag:comm}
\xy
(-10,10)*+{\X(G)}="1";
(20,10)*+{\X_\gen(G)}="2";
(-10, -10)*+{\X(G')}="3";
(20, -10)*+{\X_\gen(G')}="4";
{\ar@{->}^{\imath_G} "1"; "2"};
{\ar@{->}^{w^*} "3"; "1"};
{\ar@{->}^{(\langle w\rangle)_*} "2"; "4"};
{\ar@{->}_{\imath_{G'}} "3"; "4"};
{\ar@{=>}^{\sigma} (5,3);(5,-3)} ; 
\endxy
\end{equation}
2-commutes. 
\end{lemma}

\begin{proof}
  We will show that the functors
  $(\langle w\rangle)_* \circ \imath_G \circ w^*$ and $\imath_{G'}$
  are isomorphic by directly constructing a natural isomorphism
  $\sigma: (\langle w\rangle)_* \circ \imath_G \circ w^* \Rightarrow
  \imath_{G'}$.

  We have already seen in the proof of Lemma~\ref{lemma:key} that for
  any multiplicative vector field $u:G'\to TG'$ the bibundles
  $(\langle w \rangle)_* (\langle u|_G \rangle) = T\langle w \rangle
  \circ (\langle u|_G\rangle \circ \langle w\rangle \inv)$ and
  $\langle u \rangle$ are isomorphic.  We need to be more precise
  about choosing these isomorphisms: we have to  making sure that they are
  isomorphisms in $\X_\gen (G')$ from $(\langle w \rangle)_* (\langle
  u|_G\rangle, \alpha_{\langle
  u|_G\rangle})$ to $(\langle u \rangle , \alpha_{\langle u\rangle})$
and that they assemble into a natural transformation.
As before we may assume that $G$ is an open subgroupoid of $G'$ and
$w:G\to G'$ is an inclusion.

The category $\X_\gen (G')$ implicitly depends on the component
$\tilde{\pi}_{G'}$ of the transformation $\tilde{\pi}: T^\Bi \Rightarrow
\id_\Bi$.  Recall that the bibundle $\tilde{\pi}_{G'}$ is isomorphic to
the bibundle $\langle \pi_{G'}\rangle$.     Replacing $\tilde{\pi}_{G'}$
in the definition of $\X_\gen(G')$ by $\langle \pi_{G'}\rangle$  and
of $\tilde{\pi}_G$ by $\langle \pi_G \rangle$ in the  definition of
$\X_\gen (G)$ results in isomorphic categories.    Consequently we may
assume that $\tilde{\pi}_{G'} = \langle \pi_{G'}\rangle$ and
$\tilde{\pi}_G = \langle \pi_G\rangle$.

By definition $\langle w\rangle $ is the fiber product
$G_0\times _{w, G_0', t} G_1'$.  Here as usual $t:G_1'\to G_0'$ is the
target map.  Since the fiber product is defined by its universal
property we may assume that $\langle w\rangle = t\inv (G_0) =\{ \gamma
\in G_1'\mid t(\gamma) \in G_0\}$. Then  the left anchor $a^L_{\langle
  w\rangle}$ is the restriction of the target  map $t$ and the right
anchor  is the restrictions of the  source
map $s$.  The inverse of $\langle w\rangle$ is then $t\inv (G_0)$
with the left and right anchors reversed.

We need some notation for elements of the composites of two bibundles.
Given two composable bibundles $A\xrightarrow{P} B\xrightarrow{Q} C$
their composite is $Q\circ P = (P\times_{a^R_P, B_0, a^L_Q} Q)/B$.
So an element of $Q\circ P$ is the $B$-orbit $[p,q]$ of an element
$(p, q)$ in the fiber product $P\times_{a^R_P, B_0, a^L_Q} Q \subset
P\times Q$.  For
example given a multiplicative vector field $u:G'\to TG'$, an element
of the composite
$\langle \pi_{G'} \rangle \circ \langle u\rangle$ is of the form
\[
[(x, u(x) \xleftarrow{\mu} \dot{y}), (\dot{y}, \pi_{G'} (\dot{y})
\xleftarrow{\nu} z]
  \]
  for some $x\in G_0'$, $\mu\in TG_1$, $\dot{y}\in TG_0$, $\nu \in
  G_1'$ and $z\in G_0'$.
The natural isomorphism $\alpha_{\langle u\rangle} : \langle \pi_{G'} \rangle
\circ \langle u\rangle \to G_1' = \langle \id_{G'}\rangle$ is then
given by
\[
  \alpha_{\langle u\rangle} \left(
[(x, u(x) \xleftarrow{\mu} \dot{y}), (\dot{y}, \pi_{G'} (\dot{y})
\xleftarrow{\nu} z]
    \right) = \quad x\xleftarrow{\pi_{G'} (\mu)\circ \nu} z
  \]
For any morphism $\beta:u\Rightarrow v$ in the category $\X(G')$ the
isomorphism of bibundles $\langle \beta \rangle : \langle u\rangle \to
\langle v \rangle$ is given by
\[
  \langle \beta \rangle \left(x, u(x) \xleftarrow{\Gamma} \dot{z}
    \right) = (x, v(x) \xleftarrow{\beta(x) \circ \Gamma}\dot{z}).
  \]
  Given an element
\[
  [[y\xrightarrow{\gamma}x, (x,
 u(x)\xleftarrow{\dot{\mu}}\dot{z})], \dot{z}\xleftarrow{\dot{\nu}}
 \dot{d}]
\]
of 
\[
  \left (\langle w \rangle \inv \times _{G_0} \langle
     u|_G\rangle )/G\right) \times_{TG_0} T\langle w \rangle ) /TG =
   T\langle w \rangle \circ (\langle u|_G \rangle\circ \langle
   w\rangle \inv )\] we define
 
  \[
\sigma_u( [[y\xrightarrow{\gamma}x, (x,
u(x)\xleftarrow{\dot{\mu}}\dot{z})],
\dot{z}\xleftarrow{\dot{\nu}}
   \dot{d}] ) :=
   ( y, u(y) \xleftarrow{u(\gamma\inv)\circ \dot{\mu} \circ \dot{\nu}}
    \dot{d}) \in \langle u\rangle.
   \]
 It is not very difficult to check that $\sigma_u$ is well-defined
 isomorphism of the bibundles.
For example, the inverse of $\sigma_u$ can be constructed as follows.
Consider a point $(y, u(y)\xleftarrow{\tau} \dot{d}) \in \langle
u\rangle$.  Since $w:G\to G'$ is an essential equivalence so is $Tw:
TG\to TG'$.  In particular $Tw $ is essentially surjective.  So there
is an arrow $\dot{d} \xrightarrow{\dot{\nu}} \dot{z}$ in $TG_1'$ with
$\dot{z} \in TG_0$.    Similarly there is an arrow
$y\xrightarrow{\gamma} x$ in $G_1'$ with $x\in G_0$.  We set
\[
\sigma_u\inv (y, u(y) \xleftarrow{\tau} \dot{d}) :=
[[y\xrightarrow{\gamma} x, (x, u(x)\xleftarrow{ u(\gamma) \circ \tau
  \circ \dot{\nu}    })\dot{z} ], \dot{z} \xleftarrow{\dot{\nu}\inv} \dot{d} ]. 
  \]
One checks further that 
 \[
\alpha_{\langle u\rangle} \circ \sigma_u = \langle w\rangle_*
(\alpha_{\langle u|_G\rangle}).
\]
Hence $\sigma_u$ is a morphism in $\X_\gen (G')$ from
\[
(\langle w \rangle_* \circ \imath_G \circ w^*)(u)  = \langle w\rangle
_* (\langle u|_G \rangle, \alpha_{\langle u|_G\rangle})
\]
to $(\langle u\rangle, \alpha_{\langle u\rangle) })= \imath_{G'} (u)$.\\

It remains to check that the isomorphisms $\{\sigma_u\}$ are
components of a natural isomorphism
\[
\sigma: (\langle w \rangle_* \circ
\imath_G \circ w^*) \Rightarrow \imath_{G'}.
\]
Thus for a 
a morphism   $\beta:u\Rightarrow  v$  in $\X(G')$ we need to check
that
\[
\imath_{G'} (\beta)\circ \sigma_u = \sigma_v \circ \left( (\langle w \rangle_* \circ
\imath_G \circ w^*) (\beta)\right).
  \]
The isomorphism of bibundles 
\[
\left((\langle w \rangle_* \circ
\imath_G \circ w^*) (\beta)\right)=  \langle w \rangle_* (\langle \beta|_G\rangle): \langle w\rangle_* (
\langle w\rangle_* (\langle u|_G \rangle) \to \langle w\rangle_* (
\langle w\rangle_* (\langle v|_G \rangle)
\]
is given by
\[
[[y\xrightarrow{\gamma}x, (x,
u(x)\xleftarrow{\dot{\mu}}\dot{z})],
\dot{z}\xleftarrow{\dot{\nu}}
   \dot{d}] ) \mapsto [[y\xrightarrow{\gamma}x, (x,
v(x)\xleftarrow{\beta(x) \circ \dot{\mu}}\dot{z})],
\dot{z}\xleftarrow{\dot{\nu}}
   \dot{d}] )
 \]
 Hence
 \[
   \sigma_v \left (\langle w\rangle_* (\langle \beta|_G\rangle )) \,\left(
[[y\xrightarrow{\gamma}x, (x,
u(x)\xleftarrow{\dot{\mu}}\dot{z})],
\dot{z}\xleftarrow{\dot{\nu}}
   \dot{d}] 
 \right)\right) = (y, v(y) \xleftarrow{
 v(\gamma\inv) \circ \beta(x) \circ \dot{\mu}\circ\dot{\nu}} \dot{d})
\]
On the other hand
\[
\langle \beta\rangle (\sigma_u (
[[y\xrightarrow{\gamma}x, (x,
u(x)\xleftarrow{\dot{\mu}}\dot{z})],
\dot{z}\xleftarrow{\dot{\nu}}
\dot{d}] ) = (y, v(y) \xleftarrow{
\beta(y)\circ u(\gamma\inv) \circ \dot{\mu} \circ \dot{\nu}
  } \dot{d})
\]
Since $\beta$ is a natural transformation
\[
\beta(y)\circ u(\gamma\inv)  = v(\gamma\inv) \circ \beta(x) 
  \]
  for any arrow $y\xrightarrow{\gamma} x$ of $G'$.  Therefore
  \[
\sigma_v \circ (\langle  w\rangle _* (\langle \beta|_G \rangle)) = \langle \beta \rangle \circ \sigma_u
    \]
  for any morphism $\beta:u\Rightarrow v$ in $\X(G')$.  We conclude
  that $\{ \sigma_u\}$ are components of the desired natural isomorphism.
   \end{proof}

\section{Categories of vector fields on stacks and Lie 2-algebras} 
\label{sec:5}

We recall Hepworth's construction \cite{Hepworth} of the category
of vector fields $\Vect(\calA)$ on a stack $\calA$.  The first step is
to extend the tangent functor $T:\Man\to \Man$ on the category of
manifolds to a functor $T^\Stack: \Stack \to \Stack$ on the 2-category
of stacks over manifolds along the Yoneda embedding $y:\Man \to
\Stack$. This results in a 2-commuting diagram
\[
\xy
(-11,8)*+{\Stack}="1";
(11,8)*+{\Stack}="2";
(-11, -8)*+{\Man}="3";
(11, -8)*+{\Man}="4";
{\ar@{->}^{y} "3"; "1"};
{\ar@{->}_{y} "4"; "2"};
{\ar@{->}^{T^\Stack}"1"; "2"};
{\ar@{->}_{T}"3"; "4"};
{\ar@{=>} (-2,-2);(2,2)} ; 
\endxy 
\]
and there is a natural transformation $\pi:T^\Stack
\Rightarrow \id_\Stack$.

\begin{definition}[Hepworth]\label{def:cvfH}
  The objects of the {\sf category of vector fields} $\Vect(\calA)$ on a
  stack $\calA$ are pairs $(v,\alpha_v)$ where $v:\calA\to T^\Stack\calA$ is
  a 1-morphism of stacks and $\alpha_v: \pi_\calA\circ v \Rightarrow
 \id_\calA$ is a 2-morphism. A morphism in $\Vect(\calA)$ from $(v,\alpha_v)$ 
to $(u,\alpha_u)$ is a
  2-morphism $\beta:v\Rightarrow u$ so that
\[
\alpha_u\circ (\pi_\calA\star \beta) = \alpha_v.
\]
Here $\circ$ is the vertical composition and $\star$ is whiskering.
\end{definition}

Next recall that for any Lie groupoid $G$ there is a stack $\B G$ of
principal $G$-bundles.  The assignment
\[
G\mapsto \B G
\]
can be promoted to a functor $\B$ in different ways depending on which
source 2-category one chooses.  Hepworth takes the source to be the
2-category $\LieGpd$ of Lie groupoids, smooth functors and natural
isomorphisms and considers the functor
\begin{equation}\label{eq:B-Hep}
\B: \LieGpd \to \Stack.
\end{equation}
The essential image of this functor consists of the 2-category
$\GStack$ of geometric stacks.  The functor $\B$ is faithful but not
full.  In particular the functor $\B$ maps essential equivalences of
Lie groupoids (which need not be invertible in $\LieGpd$, even weakly) to
isomorphisms of stacks.\footnote{Recall that by tradition a weakly
  invertible 1-morphism of stacks is called an {\em isomorphism.}}  The
tangent functor $T:\Man \to \Man$ is easily extended to a functor
$T^\LieGpd:\LieGpd \to \LieGpd$.  We have a natural transformation
$\pi^\LieGpd: T^\LieGpd \Rightarrow \id_\LieGpd$.

Hepworth proves \cite[Theorem~3.11]{Hepworth} that there is a natural
isomorphism
\begin{equation} \label{eq:05.2}
\B \circ T^\LieGpd \Leftrightarrow T^\Stack \circ \B.
\end{equation}
Consequently given a vector field $v:G\to TG$ on a Lie groupoid $G$ we
get a map of stacks
\[
\B v: \B G \to \B TG.
\]
Composing $v$ with the isomorphism $\B TG \to T^\Stack (\B G)$ gives
us a functor that we again denoted by $\B v:\B G \to T^\Stack (\B G)$. 
This determines an object $(\B v, a_{Bv})$ in the category $\Vect (\B G)$ of vector fields
on the stack $\B G$.   Hepworth shows that the assignment
\[
v \mapsto (\B v, a_{Bv})
\]
can be promoted to a functor 
\begin{equation}\label{eq:functor-vf}
\X(G) \to \Vect (\B G).
\end{equation}
Here as before $\X(G)$ denotes the category of vector fields on a Lie
groupoid $G$ (see Definition~\ref{def:XG}).  He proves in
\cite[Theorem~4.15]{Hepworth} that if the groupoid $G$ is proper then
the functor \eqref{eq:functor-vf} is an equivalence of
categories.\footnote{The hypothesis that the groupoid $G$ is proper is
  not explicit in the statement of \cite[Theorem~4.15]{Hepworth}.
  However the proof depends on several lemmas: \cite[4.11, 4.12, 2.11, 2.12]{Hepworth}.
  In particular the proof uses the existence of partitions of unity
  and Weinstein-Zung linearization, both of which require properness.}

Another important consequence of the existence of the isomorphism
\eqref{eq:05.2} is that for any geometric stack $\calA$ the tangent
stack $T^\Stack \calA$ is geometric as well.\\

We can promote the assignment $G\to \B G$ to a functor out of a
different bicategory, which at a slight risk of confusion we will
again denote by $\B$. Namely we can choose as our source the
bicategory $\Bi$ of Lie groupoids, bibundles and isomorphisms of
bibundles.  The advantage is that the functor
\[
\B: \Bi \to \Stack
\]
is fully faithful: for Lie groupoids $G$ and
$H$, the functor
\[
\B : \Hom_\Bi (G, H) \to \Hom_\Stack (\B G, \B H)
\]
is an equivalence of categories.  Consequently the functor
\[
\B: \Bi \to \GStack
\]
is an equivalence of bicategories \cite{Blo}.  It is not hard to adapt
\cite[Theorem~3.11]{Hepworth} to this setting:  the diagram
\[
\xy
(-16,8)*+{\Bi}="1";
(16,8)*+{\Bi}="2";
(-16, -8)*+{\GStack}="3";
(16, -8)*+{\GStack}="4";
{\ar@{->}_{\B} "1"; "3"};
{\ar@{->}^{\B} "2"; "4"};
{\ar@{->}^{T^\Bi}"1"; "2"};
{\ar@{->}_{T^\Stack}"3"; "4"};
{\ar@{=>} (-2,-2);(2,2)} ; 
\endxy 
\]
2-commutes.  It will be convenient for us to choose a weak inverse
$\B\inv :\GStack \to \Bi$ and consider the functor
\[
T^\GStack:\GStack \to \GStack, \qquad T^\GStack:= \B \circ T^\Bi\circ \B\inv,
\]
which by construction is isomorphic to Hepworth's functor $T^\Stack$
restricted to geometric stacks.  As in the case of $T^\Stack$ we have
a transformation $\pi:T^\GStack\Rightarrow \id_\GStack$: $\pi:= \B
\star \tilde{\pi} \star \B\inv.$

Given a geometric stack $\calA$ we now define a category of vector
fields $\Vectp(\calA)$ on $\calA$ as follows (compare with
Definition~\ref{def:cvfH}).

\begin{definition}
  \label{def:cvf_on_stack}
  The {\sf category of vector fields} $\Vectp(\calA)$ on a geometric stack $\calA$,
 has as objects pairs $(X,\alpha_X)$ where
  $X:\calA \to T^\GStack \calA$ is a 1-morphism of stacks and
  $\alpha_X: \pi_\calA\circ X \Rightarrow \id_\calA$ is a
  2-morphism. A morphism from $(X,\alpha_X)$ to $(Y,\alpha_Y)$ in
  $\Vectp(\calA)$ is a 2-morphism $\beta:X\Rightarrow Y$ so that
\[
\alpha_Y\circ (\pi_\calX \star \beta) = \alpha_X.
\]
\end{definition}

It is easy to see that for a geometric stack $\calA$
the categories $\Vect(\calA)$ and $\Vectp(\calA)$ are equivalent (and even
isomorphic).  For us there are several advantages in working with
$\Vectp(\calA)$.  First of all, the functor $\Vectp$ is more explicit than
$T^\Stack$: the latter involves 2-limits and stackification.  Additionally 
the following result is easy to
prove:

\begin{lemma}\label{lem:2.12n}
  For a Lie groupoid $G$ the classifying stack functor $\B:\Bi \to
  \GStack$ induces an equivalence of categories
\[
(\B_*)_G: \X_\gen(G)\to \Vectp(\B G),
\]
where  $\X_\gen(G)$  is the category of generalized vector fields (Definition~\ref{ex:Xgen}).
\end{lemma}
\begin{proof}
  Consider a generalized vector field $(P,\alpha_P)$ on the Lie
  groupoid $G$.  By definition we have an isomorphism
  $\alpha_P:\tilde{\pi}_G \circ P \Rightarrow
  \langle\id_G\rangle$ of bibundles.  Apply the classifying stack
  functor $\B$ to the 2-morphism $\alpha_P$.  We get the 2-morphism of
  stacks
\[
\B \alpha_P:\B(\tilde{ \pi}_G
  \circ P )\Rightarrow \B\langle\id_G\rangle.
\]
Since $\B$ is functor between bicategories, we have canonical 2-arrows
$\B\langle\id_G\rangle  \Rightarrow\id_{\B G}$ and
$\B\langle \pi_G\rangle\circ \B P \Rightarrow \B(\tilde{\pi}_G
\circ P )$.  Note that these 2-morphisms are 2-isomorphisms since all
2-arrows in the 2-category of stacks are invertible.
Composing the three 2-arrows we get a 2-arrow
\[
\B\tilde{ \pi}_G\circ \B P \Rightarrow\id_{\B G}
\]
which we denote by $\alpha_{\B P}$.  By definition the pair $(\B P,
\alpha_{\B P})$ is an object of $\Vectp(\B G)$.

Similarly a morphism $\beta :(P,\alpha_P)\rightarrow (Q, \alpha_Q)$ in
$\X_\gen (G)$ gives rise to a morphism $\B \beta: \B P \Rightarrow \B
Q$.  One checks that
\[
\alpha_{\B Q} \circ (\pi_{\B G}\star \B \beta) = \alpha _{\B Q}.
\]
Consequently $\B\beta$ is a morphism in $\Vectp(\B G)$ from
$(\B P,\alpha_{\B P})$ to $(\B Q, \alpha_{\B Q})$.  We therefore get a
functor 
\[
(\B_*)_G: \X _\gen(G) \to \Vectp(\B G).
\]
A weak inverse $\B \inv :\GStack \to \Bi$ gives rise to the functor
\[
((\B\inv)_*)_G: \Vectp (\B G)\to  \X_\gen (G)
\]
in the other direction.  The induced functors $(\B_*)_G$ and $((\B \inv)_*)_G$ are
weak inverses of each other.
\end{proof}

\begin{remark}\label{rmrk:0proper}
  Suppose $G$ is a Lie groupoid.  Tracing carefully through the
  definitions: of the map $\X(G) \to \Vect(\B G)$ (this map is defined
  by Hepworth), of the equivalence $\Vectp(\B G)\to \Vect(\B G)$, of
  the %
  inclusion $\X(G)\hookrightarrow
  \X_\gen (G)$ and of the map $(\B_*)_G: \X_\gen (G)\to \Vectp (\B G)$ one
  can show that the diagram
\[
\xy
(-16,8)*+{\X(G)}="1";
(16,8)*+{\Vect(\B G)}="2";
(-16, -8)*+{\X_\gen (G)}="3";
(16, -8)*+{\Vectp (\B G)}="4";
{\ar@{^{(}->} "1"; "3"};
{\ar@{->}_{\B _*} "3"; "4"};
{\ar@{->}^{ }"1"; "2"};
{\ar@{->}^{\sim}"4"; "2"};
{\ar@{=>} (0,-2);(0,2)} ; 
\endxy 
\]
2-commutes.  By Lemma~\ref{lem:2.12n} the functor $(\B_*)_G: \X_\gen(G)\to
\Vectp(\B G)$ is  an equivalence of categories.  If the Lie groupoid $G$ is proper, the functor $\X(G) \to \Vect (\B G)$ is
an equivalence of categories by \cite[Theorem~4.15]{Hepworth}.
Consequently the functor  $\X(G) \to \X_\gen(G)$ has to be an
equivalence of categories in this case. In general Theorem~\ref{thm:04.3} tells us that the functor $\X(G) \to
\X_\gen(G)$ is an equivalence of categories for any Lie groupoid
$G$. Consequently the functor $\X(G) \to \Vect (\B G)$ is always an
equivalence of categories, regardless of whether the Lie groupoid $G$
is proper or not.  Thus Theorem~\ref{thm:04.3} generalizes
\cite[Theorem~4.15]{Hepworth}.
\end{remark}

We now address the issue of giving the category of vector fields
$\Vectp(\calA)$ on a geometric stack $\calA$ the structure of a Lie
2-algebra.  We may proceed as follows.  Choose an atlas $G_0\to \calA$
on the stack $\calA$.  The atlas induces an isomorphism of stacks
$\calA\xrightarrow{p} \B G$, where $G$ is the Lie groupoid defined by
the atlas.  The isomorphism $p$ induces an equivalence of categories
\[
p_*: \Vectp(\B G) \to \Vectp(\calA).
\]
By Lemma~\ref{lem:2.12n} the classifying stack functor induces an
equivalence of categories
\[
(\B_*)_G: \X_\gen(G) \to \Vectp(\B G).
\]
By Theorem~\ref{thm:04.3} the inclusion
\[
\imath_G :\X(G)\to \X_\gen(G)
\]
is an equivalence of categories. %
Consequently the 
composite functor $\phi_G:\X(G)\to \Vectp(\calA)$,
\[
\phi_G:= p_*\circ( \B_*)_G \circ \imath_G 
\]
is an equivalence of categories.  By Lemma~\ref{lem:03.4} the category
$\X(G)$ of multiplicative vector fields has a natural structure of a
strict Lie 2-algebra.  We may view the functor $\phi_G$ as a kind of a
``Lie 2-algebra atlas'' on the category $\Vectp(\calA)$.

What happens if we choose a different atlas $q:H_0\to \calA$ on the
stack $\calA$? By the same argument as above we get an equivalence of
categories $\phi_H :\X(H)\to \Vectp(\calA)$, which we may view as a
different ``Lie 2-algebra atlas'' on the category $\Vectp(\calA)$.  We
would like the two ``atlases'' to be compatible.  In particular we
would like to make sure that the functor
\[
\phi_G\inv \circ \phi_H:\X(H)\to \X(G)
\]
underlies a Morita equivalence of Lie 2-algebras.  At the very least we would
like the Lie 2-algebras $\X(G)$ and $\X(H)$ to be Morita
equivalent in general. That is, we would like there
to exist a 
weakly invertible 1-morphism in the bicategory $\Lietwo$ from the Lie
2-algebra $\X(H)$  to the Lie 2-algebra $\X(G)$.

To address this issue we study the functoriality of the assignment
$\calA\mapsto \X(G)$ of a Lie 2-algebra of vector fields to a
geometric stack by  a choice of an atlas $G_0\to \calA$.
Consider the 2-category $\GStack_\iso$ of geometric stacks, {\em
  isomorphisms} of stacks and 2-morphisms of stacks.  The classifying
stack functor
\[
\B :\Bi \to \GStack
\]
restricts to an equivalence of bicategories
\[
\B: \Bi_\iso \to \GStack_\iso .
\]
A choice of a weak inverse $\B\inv$ of $\B$ amounts to choosing an
atlas for each geometric stack.  Once the inverse $\B\inv$ is chosen,
we have the composite functor
\[
\GStack_\iso \xrightarrow{\B\inv} \Bi_\iso \xrightarrow{\X} \Lietwo.
\]
By construction, for a stack $\calA$ the Lie 2-algebra $\X(\B\inv(\calA))$ is the Lie 2-algebra of vector fields on the Lie groupoid $G =
\B\inv (\calA)$.  By the discussion above the category underlying the
Lie 2-algebra $\X(\B\inv (calA))$ is equivalent to the category of
vector fields $\Vectp(\calA)$ on the stack $\calA$.

A different choice of a weak inverse $(\B\inv)'$ of $\B$ amounts to
choosing a possibly different atlas for each geometric stack.  Once
$(\B\inv)'$ is chosen we have a natural isomorphism $\alpha:\B\inv
\Rightarrow (\B\inv)'$.  For each geometric stack $\calA$ the
component $\alpha_\calA$ of the natural transformation $\alpha$ is an
invertible bibundle
\[
\alpha_\calA:\B\inv (\calA) \to (\B\inv)' ( \calA).
\]
Applying the functor $\X: \Bi_\iso \to \Lietwo$ to $\alpha_\calA$ we
get an invertible bibundle
\[
\X(\alpha_\calA): \X (\B\inv ( \calA))\to \X((\B\inv)'(\calA))
\]
in the bicategory $\Lietwo$.   

One can be fairly explicit as to what the bibundle $\X(\alpha_\calA)$
actually is.  Namely let $G_0 \to \calA$ be the atlas giving rise to
the Lie groupoid $G=\B\inv (\calA)$ and $H_0 \to \calA$ be the atlas
giving rise to $H= (\B\inv)' (\calA)$.  Then the total space of the
bibundle $\alpha_\calA: G\rightarrow H =(\B\inv)' (\calA)$
represents the fiber product $G_0 \times _\calA H_0$.  The linking
groupoid $G*_{\alpha_\calA} H$ is the groupoid corresponding to the
atlas $G_0 \sqcup H_0 \to \calA$. The linking groupoid comes with two
canonical essentially surjective open embeddings
\[
i_G :G\hookrightarrow G*_{\alpha_\calA} H \quad \textrm{and}\quad 
i_H :H\hookrightarrow G*_{\alpha_\calA} H. 
\]
By Lemma~\ref{lemma:key} the
pullback/restriction functors
\[
i_G^*: \X(G*_{\alpha_\calA} H) \to \X(G), \qquad
 i_H^*: \X(G*_{\alpha_\calA} H) \to \X(G)
\]
are 1-morphisms of Lie 2-algebras that are fully faithful and
essentially surjective.  Hence the bibundle $\langle i_G^*\rangle$ is
invertible in the bicategory $\Lietwo$.  On the other hand, as was
noted in the proof of Theorem~\ref{thm:loc}, the bibundles $\langle
i_H\rangle \circ \alpha_\calA$ and $\langle i_G \rangle$ are
isomorphic.  Hence 
\[
\X (\langle
i_H\rangle) \circ \X(\alpha_\calA)  \simeq \X(\langle i_G \rangle).
\]
By construction of the functor $\X: \Bi_\iso\to \Lietwo$ we have
\[
\X(\langle i_G\rangle) = \langle i_G^*\rangle \inv \qquad \textrm{and }
\qquad 
\X(\langle i_H\rangle) = \langle i_H^*\rangle \inv. 
\]
Hence 
\[
\X(\alpha_\calA) \simeq \langle i_H^* \rangle \circ \langle i_G^* \rangle \inv.
\]

\section{Lie 2-algebras of vector fields on stacks and their 
underlying categories} 
\label{sec:5.5}

In the previous section we constructed a functor
\[
\X \circ \B\inv : \GStack_\iso\to \Lietwo.
\]
Recall that there is a forgetful functor $u:\Lietwo\to \Gpd$ that
assigns to a Lie 2-algebra its underlying groupoid.  Therefore for
every geometric stack $\calA$ we have the groupoid $(u\circ \X \circ
\B\inv)(\calA)$.  We should make sure that this groupoid is equivalent
to the groupoid of vector fields $\Vectp(\calA)$ (and hence to
Hepworth's groupoid $\Vect(\calA)$ of vector fields on the stack
$\calA$).

We start by promoting the assignment $\calA \mapsto \Vectp(\calA)$ to
a functor
\[
\Vectp:\GStack_\iso \to \Gpd
\]
whose source is the 2-category of geometric stacks and isomorphisms and whose target is the
(2,1)-category $\Gpd$ of (small) groupoids. We then prove the following theorem:

\begin{theorem}\label{thm:5.5.1}

The diagram of (2,1)-bicategories and functors
\[
\xy
(-15,8)*+{\GStack_\iso}="1";
(15,8)*+{\Gpd}="2";
(-15, -8)*+{\Bi_\iso}="3";
(15, -8)*+{\Lietwo}="4";
{\ar@{->}^{\B} "3"; "1"};
{\ar@{->}_{u} "4"; "2"};
{\ar@{->}^{\Vectp}"1"; "2"};
{\ar@{->}_{\X}"3"; "4"};
{\ar@{=>} ^{\Upsilon}(2,-2);(-2,2)} ; 
\endxy 
\]
2-commutes.   Here as above $\Gpd$ denotes the (2,1) category of groupoids, functors and natural isomorphisms, and $u:\Lietwo\to \Gpd$ denotes the functor that assigns
to each Lie 2-algebra its underlying groupoid.  The 1-components of the transformation $\Upsilon$ are weakly invertible functors (i.e., equivalences of categories).  In particular for a
geometric stack $\calA$ the category underlying the Lie 2-algebra
$(\X\circ \B\inv)\, (\calA)$
is equivalent to Hepworth's category $\Vect(\calA)$
of vector fields on the stack. 
%
%
%
%
\end{theorem}

We now construct the 2-functor $\Vectp:\GStack \to \Gpd$  (see
Appendix~\ref{sec:app} for a definition of a 2-functor).  An isomorphism 
$f:\calA_1\to \calA_2$ of stacks induces an equivalence of categories
\[
f_*:\Vectp(\calA_1) \to \Vectp (\calA_2):
\]
one adapts the proof of Lemma~\ref{lemma:equiv_of_vf2} to the setting of
geometric stacks.
Note that if $f= \id_\calA$ we may take $f_* = \id_{\Vectp(\calA)}$.  

Given isomorphisms $f:\calA_1 \to \calA_2$ and $g:\calA_2 \to \calA_3$ of
stacks we get equivalences of categories: $(g\circ f)_*$ and $g_*
\circ f_*$.  We need to produce a natural transformation $\mu_{gf}:g_*
\circ f_* \Rightarrow (g\circ f)_*$.  So given an object $(v, a_v)$ of
$\Vectp(\calA_1)$ we need to produce a 2-cell
\[
(\mu_{gf})_{(v, a_v)} : g_*(f_* (v, a_v))\Rightarrow (g\circ f)_* (v,a_v)
\]
in the category $\Vectp (\calA_3)$.  By the (adapted) proof of Lemma~\ref{lemma:equiv_of_vf2}
\[
 g_*(f_* (v)) = T^\GStack g \circ (T^\GStack f \circ v \circ f\inv )\circ g\inv.
\]
Since $T^\GStack $ is a (pseudo-) functor, there is a natural isomorphism 
\[
T^\GStack g \circ T^\GStack f \Rightarrow T^\GStack (g\circ f). 
\]
Consequently there is an isomorphism 
\[
 T^\GStack g \circ (T^\GStack f \circ v \circ f\inv )\circ g\inv \Rightarrow
T^\GStack (g\circ f) \circ v \circ (g\circ f)\inv.
\]
This isomorphism is the desired 2-cell $(\mu_{gf})_{(v, a_v)}$.  We
are now ready to describe the functor $\Vectp$.  To a geometric stack
$\calA$ it assigns the category $\Vectp(\calA)$.  To an arrow $f:\calA_1\to
\calA_2$ it assigns the equivalence of categories $\Vectp(f):= f_*$.
Additionally for each pair $(g,f)$ we have a natural isomorphism
$\mu_{gf}: g_* \circ f_* \Rightarrow (g\circ f)_*$ constructed above.
Proceeding similarly (and keeping track of the coherence data) we can
promote the assignment
\[
\Bi_\iso \ni (G\xrightarrow{P}H)\mapsto (\X_\gen(G) \xrightarrow{P_*}
\X_\gen (H) )
\]
to a functor
\[
\X_\gen : \Bi_\iso \to \Gpd.
\]
Lemma~\ref{lem:2.12n} now generalizes as follows:

\begin{lemma}\label{lemma:6.2}
The equivalences of categories 
\[
(\B_*)_G : \X_\gen (G)\to \Vectp (\B G)
\]
(one for each Lie groupoid $G$) assemble into a transformation
\[
\B_* :\X_\gen  \Rightarrow\B\circ \Vectp  .
\]
That is, the diagram
\[
\xy
(-15,8)*+{\GStack_\iso}="1";
(15,8)*+{\Gpd}="2";
(-15, -8)*+{\Bi_\iso}="3";
(15, -8)*+{\Gpd}="4";
{\ar@{->}^{\B} "3"; "1"};
{\ar@{=} "4"; "2"};
{\ar@{->}^{\Vectp}"1"; "2"};
{\ar@{->}_{\X_\gen}"3"; "4"};
{\ar@{<=} ^{\B_*}(-2,2);(2,-2)} ; 
\endxy 
\]
2-commutes. Here as before $\Gpd$ is the (2,1) category of groupoids, functors and natural isomorphisms.%
\end{lemma}

Next we prove
\begin{lemma}\label{lemma:6.3}
The diagram
\[
\xy
(-15,0)*+{\Bi_\iso}="1";
(15,8)*+{\Gpd}="2";
(15, -8)*+{\Lietwo}="4";
{\ar@{->}^{\X_\gen} "1"; "2"};
{\ar@{->}^{u} "4"; "2"};
{\ar@{->}_{\X}"1"; "4"};
{\ar@{<=}^{\imath} (1,2);(5,-2)} ; 
\endxy 
\]
2-commutes and the 1-components of $\imath$ are equivalences of categories.
\end{lemma}

\begin{proof}
  We have the underlying category functor $u_\strict:\Lietwo_\strict
  \to \Gpd$ which sends Lie 2-algebras to their underlying groupoids
  and morphisms of Lie 2-algebras to the underlying functors.  The
  functor $u_\strict$ sends essential equivalences of Lie 2-algebras
  to weakly invertible functors.  By the universal property of the
  localization $\langle \, \rangle :\Lietwo _\strict \to \Lietwo$ we
  get the underlying category functor $u:\Lietwo \to \Gpd$ with
  $u(\langle f \rangle ) $ isomorphic to $u_\strict(f)$ for every
  essential equivalence of Lie 2-algebras.  It follows that for any
  essential equivalence $f$ in $\Lietwo_\strict$ the functor
  $u(\langle f \rangle \inv)$ is a weak inverse of $u_\strict (f)$.
  We proved that for any essentially surjective open embedding $w:G\to
  G'$ of Lie groupoids the pullback functor $w^*:\X(G') \to \X(G)$ is
  an essential equivalence.  We defined $\X(w) = \langle w^* \rangle
  \inv$.  It follows that $u(\X(w))$ is a weak inverse of $u_\strict
  (w^*)$.

  By Lemma~\ref{lemma:oldkey} the diagram \eqref{diag:comm} 2-commutes for any 1-morphism $w:G\to G'$ in $\ESOE$.  Hence the diagram
\[
\xy
(-10,10)*+{\X(G)}="1";
(20,10)*+{\X_\gen(G)}="2";
(-10, -10)*+{\X(G')}="3";
(20, -10)*+{\X_\gen(G')}="4";
{\ar@{->}^{\imath_G} "1"; "2"};
{\ar@{->}_{u(\X(w))} "1"; "3"};
{\ar@{->}^{\langle w\rangle_* = \X_\gen (w)} "2"; "4"};
{\ar@{->}_{\imath_{G'}} "3"; "4"};
\endxy
\]
2-commutes as well.  It follows that the functors 
\[
u\circ \X, \X_\gen \circ \langle \, \rangle \in \Hom_W (\ESOE, \Gpd)
\]
are isomorphic (i.e., differ by a transformation whose components are equivalences of categories).  Here $\Hom_W (\ESOE, \Gpd)$ denotes the bicategory of
functors that send the collection $W$ of all 1-cells in $\ESOE$ to
weakly invertible functors.

By the definition of the functor $\X:\Bi_\iso \to\Lietwo$ its
precomposition with the localization functor $\langle \, \rangle
:\ESOE \to \Bi_\iso$ is isomorphic to $\X:\ESOE \to \Lietwo$.  It
follows that the functors $u\circ \X \circ \langle \, \rangle$ and
$\X_\gen \circ \langle \, \rangle$ are isomorphic in $\Hom_W (\ESOE,
\Gpd)$.  By the universal property of the localization
$\langle \, \rangle :\ESOE \to \Bi_\iso$, the functors $u\circ \X$ and $\X_\gen$ are
isomorphic in $\Hom (\Bi_\iso, \Gpd)$.
\end{proof}

\begin{proof}[Proof of Theorem~\ref{thm:5.5.1}]
This now follows directly from Lemmas~\ref{lemma:6.2} and~\ref{lemma:6.3}. 
\end{proof}

\section{Generalized vector fields on a Lie groupoid versus
  multiplicative vector fields}
\label{sec:6}

In this section we prove Theorem~\ref{thm:04.3}: for any Lie groupoid
$G$ the inclusion
\[
\imath_G: \X(G)\hookrightarrow \X_\gen(G), \qquad v\mapsto (\langle v\rangle,
\alpha_{\langle v\rangle }:\langle \pi_G\rangle \circ \langle v
\rangle \Rightarrow \langle id_G \rangle).
\]
of the category of multiplicative  vector fields into the category of
generalized vector fields is fully faithful and essentially
surjective.

\begin{remark}
In the case of {\em proper} Lie
groupoids Theorem~\ref{thm:04.3} follows from
\cite[Theorem~4.15]{Hepworth} --- see Remark~\ref{rmrk:0proper}.
\end{remark}

The fact that $\imath_G$ is fully faithful is an easy consequence of
Theorem~\ref{thm:folk1}.  We now address essential surjectivity.  We
first prove:

\begin{lemma}\label{lem:2.8.2vect}
  Let $V =\{V_1\toto V_0\}$ be a 2-vector space, $v_1,\ldots, v_s\in
  V_0$ a finite collection of objects and $\{v_i
  \xleftarrow{w_{ij}}v_j\}_{i,j=1}^s$ a collection of morphisms
  satisfying the cocycle conditions:
\begin{itemize}
\item $w_{ii} = 1_{v_i}$ for all $i$;
\item $w_{ji} = {w_{ij}}\inv $ for all $i,j$;
\item $w_{ij} w_{jk} = w_{ik}$ for all $i,j,k$. 
\end{itemize}
Then for any $\lambda_1,\ldots, \lambda_s \in [0,1]$ with $\sum
\lambda_k =1$ there are morphisms $v_i \xleftarrow {z_i} \sum
\lambda_k v_k$ ($i=1,\ldots, s$) with
\[
   w_{ij} = z_i z_j\inv
\]
for all $i,j$.
\end{lemma}

\begin{proof}
  By Remark~\ref{rmrk:2term-2vect_corr} the category $V$ is isomorphic
  to the action groupoid $\{U\times V_0\toto V_0\}$ where $U_0 =
  \ker(s:V_1\to V_0)$, $\partial :U\to V_0$ is $t|_U$ and the action
  of $U$ on $V_0$ is given by
\[
u\cdot v : = v+ \partial (u).
\]
Note that the multiplication/composition in 
$\{U\times V_0\toto V_0\}$ is given by
\[
(u',v+\partial(u))(u, v) = (u'+u, v)
\]
for all $v\in V_0$, $u,u' \in U$.
Consequently 
\[
(u,v)\inv = (-u, v+\partial(u)).
\]
The isomorphism $f:V\to \{U\times V_0\toto
V_0\}$ is given on morphisms by
\[
f(w) = (w-1_{s(w)}, s(w))\in U\times V_0 \qquad \textrm{for all } w\in V_1.
\]
The isomorphism $f$ followed by the projection onto $U$ sends the
morphisms $w_{ij}$ to vectors $u_{ij}\in U$.  It is easy to see that the
cocycle conditions translate into:
\begin{itemize}
\item $u_{ii} = 0$ for all $i$;
\item $u_{ji} = -u_{ij}$ for all $i,j$;
\item $u_{ik} -  u_{jk} = u_{ij}$ for all $i,j,k$. 
\end{itemize}
Moreover 
\[
\partial(u_{ij}) = v_i - v_j \qquad \textrm{for all } i,j.
\]
Now consider 
\[
y_i = (\sum \lambda_k u_{ik}, \sum \lambda_k v_k) \in U\times V_0 
\]
and set 
\[
z_i := f\inv (y_i) \in V_1.
\]
We now verify that the $z_i$'s are the desired morphisms.
By definition the source of $y_i$ is $\sum \lambda_k v_k$.  The target of $y_i$ is 
\begin{eqnarray*}
  \partial(\sum_k \lambda_k u_{ik}) +  \sum_k \lambda_k v_k = 
  \sum_k \lambda_k \partial(u_{ik}) +  \sum_k \lambda_k v_k \\
  = \sum_k  \lambda_k (v_i - v_k) +  \sum_k \lambda_k v_k = 
\sum_k\lambda_k v_i = v_i.
\end{eqnarray*}
Hence $z_i$ is an arrow from $\sum \lambda_k v_k$ to $v_i$. Finally
\begin{eqnarray*}
  y_i y_j\inv = 
(\sum_k \lambda_k u_{ik}, \sum \lambda_k v_k)(- \sum_k \lambda_k u_{jk}, v_j) \\
= (\sum_k \lambda_k (u_{ik} - u_{jk}), v_j) = (\sum_k \lambda_k u_{ij}, v_j)=
(u_{ij}, v_j),
\end{eqnarray*}
and so $z_iz_j\inv = w_{ij}$ as desired.
\end{proof}

\begin{proposition} \label{thm:1} Let $G=\{G_1\toto G_0\}$ be a Lie
  groupoid, $U_0\subset G_0$ an open submanifold and $U= \{U_1\toto
  U_0\}$ the restriction of $G$ to $U_0$ (that is, $U_1$ consists of
  arrows of $G$ with source and target in $U_0$).  Given a functor
  $X:U\to TG$ together with a natural isomorphism $\alpha:
  (i:U\hookrightarrow G)\Rightarrow \pi_G \circ X $ there exists a
  functor $Y:U\to TU$ so that $\pi_U\circ Y =\id_U$ and a natural
  isomorphisms $\beta: Ti \circ Y \Rightarrow X $.
\end{proposition}
\begin{proof}
By definition of $\alpha$ the diagram
\[
\xy
(-10,10)*+{U_0}="1";
(10,10)*+{TG_0}="2";
(-10, -10)*+{G_1}="3";
(10, -10)*+{G_0}="4";
{\ar@{->}^X "1"; "2"};
{\ar@{->}_\alpha "1"; "3"};
{\ar@{->}^{\pi_G} "2"; "4"};
{\ar@{->}_t "3"; "4"};
\endxy
\]
commutes.  Hence there is a smooth map 
\[
(\alpha, X):U_0 \to G_1\times_{t, G_0, \pi}TG_0 = t^*TG_0. 
\]
Since the target map $t:G_1\to
G_0$ is a submersion, its differential 
\[
Tt_\gamma:T_\gamma G_1 \to T_{t(\gamma)}G_0 
\]
is a surjective linear map for each $\gamma\in G_1$.  Consequently the map 
\[
\Phi:TG_1 \to t^*TG_0, \qquad \Phi(\gamma, v) = (\gamma, Tt_\gamma v)
\]
is a surjective map of vector bundles over $G_1$.  Choose a smooth
section $\sigma: t^*TG_0 \to TG_1$ of $Tt$ and consider the composite
\[
\beta:= \sigma\circ (\alpha, X): U_0 \to TG_1.
\]
By construction of $\beta$  
\[
\beta(x) \in T_{\alpha(x)}G_1 \qquad \textrm{and}\qquad 
Tt_{\alpha(x)}  \beta (x) = X(x)
\]
 for any $x\in U_0$. 
We now define a functor $Y:U\to TU$.   On objects we set 
\[
Y(x) = Ts (\beta(x)).
\]
For an arrow $x\xrightarrow{\,\gamma\,}y\in U_1$ we set 
\[
Y(\gamma) = \beta(y)\inv X(\gamma)\beta(x).
\]
It is easy to check that $Y$ is indeed functor, $\beta:Ti\circ
Y\Rightarrow X$ is a natural transformation and $\pi_U \circ Y =\id_U$.
\end{proof}

\begin{proposition}\label{thm:4.6}
 Let $G$ be a Lie groupoid and 
$$
\xy
(-19,8)*+{G_1}="1";
(-19,-8)*+{G_0}="2";
(19, 8)*+{TG_1}="3";    
(19, -8)*+{TG_0}="4";
(0,5)*+{P} = "5";
{\ar@{->} (-20,6); (-20, -5)};
{\ar@{->} (-18,6); (-18, -5)};
{\ar@{->} (20,6); (20, -5)};
{\ar@{->} (18,6); (18, -5)};
{\ar@{->}_{a_P^L}  "5";"2"};
{\ar@{->}^{a_P^R}  "5";"4"};
\endxy
$$ 
be a bibundle from $G$ to the tangent groupoid $TG$ such that the
composite $\langle \pi \rangle \circ P $ is isomorphic to $\langle
id_G \rangle$ by way of a bibundle isomorphism 
\[
{\mathbf a}:\langle
\pi \rangle \circ P \Rightarrow \langle\id_G
\rangle.
\] 
  Then the left anchor $a^L_P:P\to G_0$ has a global section
  $\tau: G_0\to P$.  Moreover we may choose $\tau$ so that the
  corresponding functor $X_\tau:G\to TG$ is a multiplicative vector
  field (i.e., $\pi_G \circ X_\tau =\id_G$). Consequently the functor $\iota_G:\X(G)\to \X_\gen (G)$ of
Theorem~\ref{def:incl_functor} is essentially surjective.
\end{proposition}
\begin{proof}
  Since $a^L_P:P\to G_0$ is a surjective submersion, it has local
  sections.  Choose a collection of local sections
  $\{\sigma_i:U_0^{(i)}\to P\}$ of $a^L_P$ so that $\{U_0^{(i)}\}$ is
  an open cover of $G_0$.  It is no loss of generality to assume that
  the cover is locally finite.  Denote the restriction of the groupoid
  $G$ to $U_0^{(i)}$ by $U^{(i)}$.  That is, the manifold of objects
  of $U^{(i)}$ is $U^{(i)}_0$ and the manifold of arrows $U^{(i)}_1$
  consists of all arrows of $G$ with source and target in $U^{(i)}_0$,
  so $U_1^{(i)}:=s^{-1}(U_0^{(i)})\bigcap t^{-1}(U_0^{(i)})$.

  For each section $\sigma_i$ we get a functor $X_i: U^{(i)} \to TG$
  whose value on objects is
\[
X_i (x) = a^R_P (\sigma_i(x)).
\]
The value of $X_i$ on an arrow $y\xleftarrow{\gamma} x \in U^{(i)}_1$
is uniquely defined by the equation
\[
\gamma \cdot \sigma_i(x) = \sigma_i(y) \cdot X_i (\gamma)
\]
(see Lemma~\ref{lem:2.5n}).
We next observe that the isomorphism ${\mathbf a}:\langle \pi_G
\rangle \circ P \to \langle\id_G \rangle $ gives rise to natural
isomorphisms $\alpha_j:\pi_G \circ X_j \Rightarrow (\imath_j: U^{(j)}
\hookrightarrow G)$ where $\imath_j: U^{(j)} \hookrightarrow G$ is the
inclusion functor. This can be seen as follows.

Recall that the composite $Q\circ P$ of bibundles
$P:K\to L$ and $Q:L\to M$ is the quotient of
the fiber product $P\times_{a^R_P, L_0, a^L_Q} Q$ by the action of
$L$.  We denote by $[p,q]$ the orbit of $(p,q)\in P\times_{a^R_P, L_0,
  a^L_Q} Q$ in $Q\circ P = (P\times_{a^R_P, L_0, a^L_Q} Q)/L$.  The
bibundle $\langle \pi_G \rangle$ is the fiber product $TG_0 \times
_{\pi_G, G_0, t}G_1$ with the anchor maps $a^R_{\langle \pi_G
  \rangle}(v, \gamma) = v$, $a^L_{\langle \pi_G \rangle}(v, \gamma)=
s(\gamma)$.   Consequently in our case 
\[
\langle \pi_G \rangle \circ P 
= (P\times _{a^R_P, TG_0,a^L_{\langle \pi_G \rangle}} (TG_0 \times_{\pi_G, G_0, s} G_1))/TG. 
\]
It is convenient to identify $P\times _{a^R_P, TG_0,a^L_{\langle \pi_G
    \rangle}}(TG_0 \times_{\pi_G, G_0, s} G_1)$ with $P\times
  _{\pi_G \circ a^R_P, G_0,s} G_1$ by way of the $TG$-equivariant
  isomorphism
\[
(p, (\pi_G \circ a^R_P)(p), \gamma)\mapsto (p, \gamma).
\]
We then have a $G\times G$ equivariant diffeomorphism
\[ {\mathbf a}: (P\times _{\pi_G \circ a^R_P, G_0,s} G_1)/TG \to
G_1,\qquad [p, \gamma] \mapsto {\mathbf a} ([p,\gamma])
\]
with
\[
s({\mathbf a} ([p,\gamma])) = s(\gamma)\quad \textrm{and} \quad
t({\mathbf a} ([p,\gamma])) = a^P_L(p).
\]
A local section $\sigma_i:U_0^{(i)}\to P$ also defines a local section 
\[
\bar{\sigma}_i:U_0^{(i)}\to (P\times _{G_0} G_1)/TG 
\]
of $a^L_{\langle \pi_G \rangle \circ P}:(P\times _{G_0} G_1)/TG \to
G_0$.  It is given by
\[
\bar{\sigma}_i(x) = [\sigma_i(x), 1_{(\pi_G \circ a^R_P \circ \sigma_i) (x)}] (= 
[\sigma_i(x), 1_{\pi_G \circ X_i \,(x)}]).
\]
The arrow ${\mathbf a}(\bar{\sigma}_i (x)) \in G_1 = \langle\id_G
\rangle$ is an arrow with the target $a^L_P (\sigma_i(x)) = x$ and the
source $s (1_{\pi_G \circ X_i \,(x)}) = \pi_G \circ X_i \,(x)$.  We
define the desired natural isomorphism $\alpha_i$ by setting
\[
\alpha_i (x) = \left( {\mathbf a}(\bar{\sigma}_i (x))\right)\inv.
\]
By Proposition~\ref{thm:1} there are smooth maps $\beta_i: U_0^{(i)}
\to TG_1$ so that
\[
\pi_G \circ \beta_i  = \alpha_i 
\]
and
\[
Tt\circ \beta_i = X_i.
\]
Moreover the functors   $Y_i : U^{(i)}\to TG$ given  by 
\[
Y_i = Ts \circ \beta_i
\]
define multiplicative vector fields on each groupoid $U^{(i)}$. This
is because their images land in $TU^{(i)} \subset TG$.  In particular
$\pi_G (Y_i(x)) = x$ for all $x\in U^{(i)}_0$.

Define the local sections $\nu_i: U_0^{(i)}\to P$ of $a^L_P$ by
\[
\nu_i (x) := \sigma_i (x) \cdot \beta_i (x) 
\]
for all $x\in U_0^{(i)}$.  Then by definition
\[
a^R_P (\nu_i (x)) = Y_i (x) 
\] 
and
\[
\gamma \cdot \nu_i (x) = \nu_i (y) \cdot Y_i (x)
\]
for all arrows $y\xleftarrow{\gamma}x $.
For all $i$ and all $x\in U^{(i)}_0$ 
\begin{eqnarray*}
{\mathbf a} ([\nu_i (x), 1_{\pi_G\circ a^R_P \circ \nu_i(x))}] )=
  {\mathbf a} ([\sigma_i (x)\beta_i (x), 1_{\pi_G \circ Y_i
  (x)}])={\mathbf a} ([\sigma_i (x), \pi_G (\beta(x))]\\
={\mathbf a} ([\sigma_i (x), 1_{\pi_G \circ X_i (x)}]) \pi_G
  (\beta(x))= {\mathbf a} (\bar{\sigma}_i (x))
  \alpha_i (x) = 1_x.
\end{eqnarray*}
Hence 
\[
{\mathbf a} ([\nu_i (x), 1_{\pi_G\circ Y_i(x)}] = 1_x.
\]
Finally, we construct a global section $\tau:G_0 \to P$ of $a^L_P$ and
the corresponding global multiplicative vector fields $X_\tau:G\to TG$
using a partition of unity argument.  Choose a partition of unity
$\{\lambda_i\}$ on $G_0$ subordinate to the cover $\{U_0^{(i)}\}$.
Since the cover is locally finite it is no loss of generality to
assume that the cover is in fact finite.

Consider a point $x\in U^{(i)}_0 \cap U^{(j)}_0$.  Then 
\[
\pi_G \circ a^R_P \circ \nu_i (x) = x = \pi_G \circ a^R_P \circ \nu_j (x).
\]
Moreover
\[ 
{\mathbf a} ([\nu_i (x), 1_x]) ={\mathbf a} ([\nu_i (x), 
1_{\pi_G\circ a_R \circ \nu_i(x))}]) = 1_x.
\]
Similarly
\[ 
{\mathbf a} ([\nu_j (x), 1_x]) = 1_x.
\]
Since ${\mathbf a}$ is a diffeomorphism it follows that 
\[
[\nu_j (x), 1_x]=[\nu_i (x), 1_x]
\]
in the orbit space $(P\times _{G_0}G_1)/TG$.
Therefore there is an arrow $w_{ij}\in TG_1$ so that
\[
(\nu_i (x)w_{ij}(x) , 1_x)=(\nu_j (x), \pi_G (w_{ij}(x))1_x).
\]
Consequently
\[
\nu_i (x)w_{ij}(x) =\nu_j (x)\qquad \textrm{and}\qquad 
\pi_G (w_{ij}(x)) = 1_x,
\]
that is, $w_{ij}(x)\in T_{1_x}G_1$.  Moreover since $a^L_P:P\to G_0$
is a principal $TG_1$ bundle, the arrow $w_{ij}(x)$ with this property
is unique and depends smoothly on $x$.  Note that the source of
$w_{ij}$ is $Y_j(x)$ and the target is $Y_i(x)$.  The uniqueness of the
$w_{ij}(x)$'s implies that the collection $\{w_{ij}(x)\}$ satisfies
the cocycle conditions of Lemma~\ref{lem:2.8.2vect}.  Therefore there
exist arrows $Y_i(x) \xleftarrow{z_i(x)} \sum_k \lambda_k Y_k (x)$ with
  $z_i(x)z_j (x) \inv = w_{ij}(x)$.  A quick look at the proof of Lemma~\ref{lem:2.8.2vect}
 should convince the reader that $z_i(x)$'s depend smoothly on $x$.  

For $x\in U_0^{(i)}$ we set $\tau(x) = \nu_i(x)\cdot z_i (x)$.   Note
that for $x\in U_0^{(i)} \cap U_0^{(j)}$
\[
\nu_j(x) = \nu_i (x)w_{ij} (x) = \nu_i (x) z_i(x) z_j(x)\inv.
\]
Therefore
\[
\nu_j(x)\cdot z_j (x) = \nu_i(x)\cdot z_i (x).
\]
It follows that $\tau$ is a globally defined section of $a^L_P:P\to
G_0$. It remains to show that the corresponding functor $X_\tau:G\to
TG$ is a multiplicative vector field. By construction for each index $i$ we have a natural isomorphism $z_i:Y_i \Rightarrow
X_\tau |_{U^{(i)}}$. 
Since $z_i(x) \in T_{1_x }G_1$ and since $Y_i$
is a multiplicative vector field, the restriction $X_\tau |_{U^{(i)}}$
is also a multiplicative vector field.  We conclude that $X_\tau$ is a
multiplicative vector field globally.

We now argue that $\imath_G: \X(G) \to \X_\gen(G)$ is essentially
surjective.    Given a generalized vector fields $(P,\alpha_P)$ the
isomorphism $\alpha_P: \tilde{\pi}_G \circ P \to \langle \id_G\rangle$
defines an isomorphism $\mathbf{a}: \langle \pi_G \rangle \circ P \to
\langle \id_G\rangle$ (remember that we fixed a 2-cell $\tilde{\pi}_G
\to \langle \pi_G \rangle$ for every Lie groupoid $G$.   By the above
argument we have a multiplicative vector field $X:G\to TG$ and an
isomorphism of bibundles $\gamma : P\to \langle X \rangle$.  It
remains to show that the vector fields $(P, \alpha_P)$ and $(\langle X
\rangle, \alpha _{\langle X \rangle}$ are isomorphic.   This is not
entirely obvious since $\alpha _{\langle X \rangle } \circ
(\tilde{\pi}_G \star \gamma)$ need not equal to $\alpha_P$.  None the
less $\gamma $ can be modified to a new isomorphism $\beta:P \to
\langle X \rangle$ so that $\alpha _{\langle X \rangle } \circ
(\tilde{\pi}_G \star \beta) = \alpha _{\langle X \rangle}$.  This
follows from Lemma~\ref{lem:7.5} below. 
\end{proof}

\begin{lemma} \label{lem:7.5}  Let $G$ be a Lie groupoid, $P:G\to TG$
  a bibundle, $\alpha,\alpha' : \tilde{P}\circ P\to \langle \id_G
  \rangle$ two isomorphisms of bibundles.  Then there exists an
  isomorphism $\beta:P\to P$ of bibundles so that
  \[
 \alpha \circ (\tilde{\pi}_G \star \beta )  = \alpha'.
    \]
Moreover we may take $\beta = P\star (\alpha' \circ \alpha\inv).$
 \end{lemma}

 \begin{proof}
For any two isomorphisms $\gamma, \delta: \langle \id_G \rangle \to
\langle \id_G\rangle$ of bibundles
\[
\gamma \circ \delta = \gamma \star \delta.
 \] 
Consequently for any $\gamma: \langle \id_G \rangle \to
\langle \id_G\rangle$ we have
\begin{align*}
  (\alpha \circ (\alpha')\inv )\circ \gamma & =
  & (\alpha \circ (\alpha')\inv )\star \gamma\\
&=&(\alpha \circ 1_{\tilde{\pi}_G \circ P} \circ (\alpha')\inv ) \star
    (1_{\langle \id_G \rangle}\circ \gamma \circ 1_{\langle \id_G
    \rangle})\\
  &=& (\alpha \star 1_{\langle \id_G \rangle}) \circ ((\tilde{\pi_G}
      \circ P) \star \gamma) \circ ((\alpha')\inv \star 1 _{\langle
      \id_G \rangle})\\
  &= & \alpha \circ ( \tilde{\pi}_G \star (P\star \gamma)) \circ (\alpha')\inv.
\end{align*}
Hence if $\gamma  = \alpha ' \circ \alpha\inv $
\[
1_{\langle \id_G \rangle } =\alpha \circ (\tilde{\pi}_G \star
                               (P\star (\alpha' \circ \alpha \inv )))
                               \circ (\alpha')\inv .
\]
and therefore
\[
\alpha' = \alpha \circ (\tilde{\pi}_G \star
                               (P\star (\alpha' \circ \alpha \inv ))).
  \]
 \end{proof}  
\appendix

\section{Bicategories, functors and natural
  transformations}\label{sec:app}
In this section for the reader's convenience we record the definitions of
bicategories, (2,1)-bicategories, (pseudo-)functors,
(pseudo-)natural transformations, functors and functor bicategories.
The original reference is \cite{Benabou}.  Our presentation
closely follows \cite{Leinster}.  

\begin{definition}
A \emph{bicategory} $\calB$ consists of the following data subject
  to the following axioms:\\[2pt]
\noindent {\sc data}
\begin{itemize}
\item A collection $\calB_0$ of 0-cells (or objects).

\item For every pair $A, B$ of 0-cells a {\em category} $\Hom_\calB
  (A,B)$ of morphisms from $A$ to $B$.  The objects of $\Hom_\calB
  (A,B)$ are called 1-cells (or 1-morphisms) and are written $f:A\to
  B$ (small Latin letters). The morphisms $\Hom_\calB (A,B)$ are
  called 2-cells (or 2-morphisms) and are written $\alpha:f\Rightarrow
  g$ (Greek letters).  Note that for every 1-cell $f$ we have the
  2-cell $\id_f:f\Rightarrow f$.

  We refer to the composition of 2-cells in $\Hom_\calB (A,B)$ as a
  vertical composition and write it as $\circ$ or as blank: $\beta
  \circ \alpha \equiv \beta \alpha$.

\item For every triple of 0-cells $A,B,C$ a {\em composition} functor
\begin{eqnarray*}
  c_{ABC} = c: \Hom_\calB (B,C)\times \Hom_\calB (A,B)&\to& \Hom_\calB
  (A,C)\\
  c_{ABC} (g,f) &= &g\circ f \equiv gf, 
\qquad \textrm{ for all 1-cells } g,f \\ 
  c_{ABC}(\beta,
  \alpha)& =& \beta \star \alpha. \qquad \textrm{ for all 2-cells } \beta, \alpha
\end{eqnarray*}
We refer to the composition $\star$ of 2-cells as the {\sf horizontal
  composition}.
\item For every object $A$ of $\calB$ a 1-cell $1_A \in \Hom_\calB (A,A)$.
\item Natural isomorphisms $a_{ABCD}:c_{ABD}\circ
  (c_{BCD}\times \id)\Rightarrow c_{ACD}\circ (\id \times
  c_{ABC})$ called {\sf associators} for every quadruple of
  objects $A,B, C, D$:
\[
\xy
(-40,10)*+{\Hom_\calB (C,D) \times\Hom_\calB (B,C)\times \Hom_\calB (A,B)}="1";
(40,10)*+{\Hom_\calB (C,D)\times\Hom_\calB (A,C) }="2";
(-40, -6)*+{\Hom_\calB (B,D)\times \Hom_\calB (A,B) }="3";
(40, -6)*+{\Hom_\calB (A,D) }="4";
{\ar@{->}^{\id_{}\times c _{ABC}} "1"; "2"};
{\ar@{->}_{c_{BCD}\times \id} "1"; "3"};
{\ar@{->}^{c_{ACD}} "2"; "4"};
{\ar@{->}_{c_{ABD}} "3"; "4"}; 
{\ar@{=>}_{a_{ABCD}} (-3,0);(3,4)} ;
\endxy
\]
and in particular invertible 2-cells
\[
a_{hgf}: (hg)f \Rightarrow h(gf)
\]
for every triple of 1-cells $(h,g,f)\in \Hom_\calB (C,D)
\times\Hom_\calB (B,C)\times \Hom_\calB (A,B)$.
\item Natural isomorphisms $r_{AB}: c_{AAB}\circ (\id
  \times 1_A) \Rightarrow \id$ and $l_{AB}: c_{ABB}\circ (1_B
  \times \id) \Rightarrow \id$ called right and left {\sf unitors} for
  every pair of objects $A,B$ of $\calB$:
\[
\xy
(-20,10)*+{\Hom_\calB (A,B)}="1";
(-20, -10)*+{\Hom_\calB (A,B)\times \Hom_\calB (A,A)}="3";
(20, -10)*+{\Hom_\calB (A,B)}="4";
{\ar@{->}_{\id\times 1_A} "1"; "3"};
{\ar@{->}^{\id} "1"; "4"};
{\ar@{->}_{\qquad \quad \circ_{AAB}} "3"; "4"};
{\ar@{=>}^{r_{AB}} (-6,-6);(-2,-2)} ;  
\endxy
\qquad
\xy
(-18,10)*+{\Hom_\calB (A,B)}="1";
(-18, -10)*+{\Hom_\calB (B,B)\times \Hom_\calB (A,B)}="3";
(21, -10)*+{\Hom_\calB (A,B)}="4";
{\ar@{->}_{1_B \times \id} "1"; "3"};
{\ar@{->}^{\id} "1"; "4"};
{\ar@{->}_{\qquad \quad \circ _{ABB}} "3"; "4"};
{\ar@{=>}^{l_{AB}} (-6,-6);(-2,-2)} ;  
\endxy
\]
where $\id:\Hom_\calB (A,B)\to\Hom_\calB (A,B)$ denotes the identity functor. By abuse of notation $1_A: \Hom_\calB(A,B) \to \Hom_\calB(A,A)$
denotes the functor that takes every 2-cell to the identity 2-cell
$\id_{1_A}: 1_A \Rightarrow 1_A$.  The functor $1_B$ is defined
similarly.  Thus for every 1-cell $f\in \Hom_\calB(A,B)$ we have {\em invertible} 
2-cells
\[
r_f: f\circ 1_A \Rightarrow f\quad \textrm{ and } \quad 
l_f: 1_B \circ f \Rightarrow f.
\]
\end{itemize}
\mbox{}\\[2pt]
\noindent {\sc conditions on the data}\\[2pt]
\begin{itemize}
\item (Triangle identity) For any pair of composable 1-cells
  $C\xleftarrow{g}B\xleftarrow{f}A$ the diagram of 2-cells
\[
\xy
(-15,20)*+{(g1_B)f}="1";
(15,20)*+{(g(1_B f)}="2";
(0,0)*+{gf} ="3";
{\ar@{=>}^{a} "1"; "2"};
{\ar@{=>}_{r_g\star 1_f} "1"; "3"};
{\ar@{=>}^{\id_g\star l_f} "2"; "3"};
\endxy
\]
commutes.
\item (Pentagon identity) For any quadruple  of composable 1-cells  $E\xleftarrow{k}D\xleftarrow{h} C\xleftarrow{g}B\xleftarrow{f}A$ the diagram 
\[
\xy
(-15,20)*+{((kh)g)f}="1";
(15,20)*+{(k(hg))f}="2";
(22,5)*+{k((hg)f)}="3";
(0,-10)*+{k(h(gf))}="4";
(-22,5)*+{(kh)(gf)}="5";
{\ar@{=>}^{a\star \id_F} "1"; "2"};
{\ar@{=>}^{a} "2"; "3"};
{\ar@{=>}^{a} "3"; "4"};
{\ar@{=>}^{\id_k \star a} "1"; "5"};
{\ar@{=>}^{a} "5"; "4"};
\endxy
\]
commutes.
\end{itemize}
\end{definition}
\begin{remark}
\begin{itemize}
\item If the natural isomorphisms $a,r$ and $l$ of a bicategory
  $\calB$ are identities so that $(hg)f = h(gf)$ and $1_B\circ f = f =
  f\circ 1_A$ for all 1-cells $h,g,f:A\to B$, and similarly for the horizontal
  composition of 2-cells, then the bicategory $\calB$ is a (strict)
  {\sf 2-category}.
\item If all the 2-cells in a bicategory $\calB$ are invertible, that
  is, for any 2-cell $\alpha:f\Rightarrow g$ there is
  $\beta:g\Rightarrow f$ with $\beta \alpha = \id_f$ and $\alpha \beta
  = \id_g$ then $\calB$ is a {\sf (2,1)-bicategory}.
\end{itemize}
\end{remark}
The following definitions are special cases of more general
definitions that are specialized to (2,1)-bicategories.
\begin{definition} \label{def:App3}
A (pseudo-){\sf functor} from a (2,1)-bicategory
  $\calB$ to a (2,1)-bicategory $\calB'$ consists of the following data:
\begin{itemize}
\item A function $F = F_0:\calB_0\to \calB_0'$ on objects.
\item For any pair of objects $A,B\in \calB_0$ a functor
\[
F_{AB}:\Hom_\calB (A,B)\to \Hom_{\calB'} (FA,FB)
\]
\item  For every triple of objects $A,B,C$ of $\calB$ a natural isomorphism
\[
\mu_{ABC}: (c'_{FA,FB,FC}) \circ (F_{BC}\times F_{AB}) \Rightarrow
F_{AC} \circ (c_{ABC}):
\]
\[
\xy
(-40,10)*+{\Hom_\calB (B,C)\times \Hom_\calB (A,B)}="1";
(40,10)*+{\Hom_\calB (A,C) }="2";
(-40, -6)*+{\Hom_{\calB'} (FB,FC)\times \Hom_{\calB'} (FA,FB) }="3";
(40, -6)*+{\Hom_{\calB'} (FA,FC) }="4";
{\ar@{->}^{c _{ABC}} "1"; "2"};
{\ar@{->}_{F_{BC}\times F_{AB}} "1"; "3"};
{\ar@{->}^{F_{AC}} "2"; "4"};
{\ar@{->}_{c'_{FA,FB, FC}} "3"; "4"}; 
{\ar@{=>}_{\mu_{ABC}} (-3,0);(3,4)} ;
\endxy
\]
thus invertible 2-cells $\mu_{gf}:Fg\circ Ff \Rightarrow F(g\circ f)$
for every pair of composable 1-cells $C\xleftarrow{g} B\xleftarrow{f}
A$.  Here $c'_{FA,FB,FC}$ is the composition functor in the target
bicategory $\calB'$.
\item 2-cells $\mu_A: 1_{FA} \Rightarrow F(1_A)$ for every object $A$ of $\calB$.
\end{itemize}
\mbox{}\\[2pt]
\noindent {\sc conditions on the data}\\[2pt]
For any triple $ D\xleftarrow{h} C \xleftarrow{g} B \xleftarrow{f} A$ of
composable 1-cells the following diagrams of 2-cells commute:
\[
\xy
(-40,10)*+{(Fh\circ Fg)\circ Ff}="1";
(0,10)*+{F(h\circ g)\circ Ff}="2"; 
(40, 10)*+{F((h\circ g)\circ f) }="3";
(-40, -10)*+{Fh\circ (Fg\circ Ff) }="4";
(0, -10)*+{Fh\circ (F(g\circ f) }="5";
(40, -10)*+{F(h\circ (g\circ f)) }="6";
{\ar@{=>}^{\mu \star 1_{Ff}} "1"; "2"};
{\ar@{=>}_{\mu} "2"; "3"};
{\ar@{=>}^{a_{\calB'}} "1"; "4"};
{\ar@{=>}_{F(a_\calB)} "3"; "6"}; 
{\ar@{=>}_{1_{Fh}\star \mu} "4";"5"} ;
{\ar@{=>}_{\mu} "5";"6"} ;
\endxy
\]
(here $a_{\calB'}$ denotes the associator in the bicategory $\calB'$),
\[
\xy
(-40,10)*+{Ff\circ 1_{FA}}="1";
(0,10)*+{Ff\circ F(1_{A})}="2"; 
(40, 10)*+{F(f\circ 1_A) }="3";
(0, -15)*+{Ff }="4";
{\ar@{=>}^{1_{Ff} \star \mu} "1"; "2"};
{\ar@{=>}_{\mu} "2"; "3"};
{\ar@{=>}^{r_{\calB'}} "1"; "4"};
{\ar@{=>}_{F(r_\calB)} "3"; "4"}; 
\endxy
\]
and
\[
\xy
(-40,10)*+{1_{FB}\circ Ff}="1";
(0,10)*+{F(1_B)\circ F(f)}="2"; 
(40, 10)*+{F(1_B \circ f) }="3";
(0, -15)*+{Ff }="4";
{\ar@{=>}^{1_{\mu\star 1_{Ff}}} "1"; "2"};
{\ar@{=>}_{\mu} "2"; "3"};
{\ar@{=>}^{l_{\calB'}} "1"; "4"};
{\ar@{=>}_{F(l_\calB)} "3"; "4"}; 
\endxy
\]
Here $r_\calB$ and $r_{\calB'}$ denote the right unitors in $\calB$ and $\calB'$ respectively.   The meaning of $l_{\calB'}$ and $l_\calB$ is similar.
\end{definition}

\begin{definition} \label{def:App_nat}  
  Let $(F,\mu)$, $(G,\tau)$ be (pseudo-) functors from a
  (2,1)-bicategory $\calB$ to a (2,1)-bicategory $\calB'$.  A {\sf
  (pseudo-)  natural transformation} $\sigma: (F,\mu) \Rightarrow (G,\tau)$
  consists of the following data:
\begin{itemize}
\item 1-cells $\sigma_A:FA\to GA$ for every object $A$ of $\calB$;
\item natural isomorphisms
\[
\xy
(-20,10)*+{\Hom_\calB (A,B)}="1";
(20,10)*+{\Hom_{\calB'} (FA,FB) }="2";
(-20, -6)*+{\Hom_{\calB'} (GA,GB) }="3";
(20, -6)*+{\Hom_{\calB' }(FA,GB) }="4";
{\ar@{->}^{F_{AB}} "1"; "2"};
{\ar@{->}_{G_{AB}} "1"; "3"};
{\ar@{->}^{\sigma_B \circ \,-} "2"; "4"};
{\ar@{->}_{-\,\circ \sigma_A} "3"; "4"}; 
{\ar@{=>}_{\sigma_{AB}} (-3,0);(3,4)} ;
\endxy
\]
and thus 2-cells $\sigma_f: Gf\circ \sigma_A \Rightarrow \sigma_B \circ Ff$.
\end{itemize}
The data are subject to the following compatibility conditions --- for
any pair of composable 1-cells 
\[
C\xleftarrow{g}B\xleftarrow{f}A 
\]
the diagrams below commute:
\begin{equation}\label{eq:App1}
\xy
(-20,30)*+{Gg \circ (Gf\circ \sigma_A)}="1";
(20,30)*+{Gg\circ (\sigma_B \circ Ff) }="2";
(-30, 10)*+{(Gg\circ Gf)\circ \sigma_A }="8";
(30, 10)*+{(Gg\circ \sigma_B ) \circ Ff }="3";
(-30, -10)*+{G(g\circ f)\circ \sigma_A }="7";
(30, -10)*+{(\sigma_C \circ Fg)\circ Ff }="4";
(-20,-30)*+{\sigma_C \circ F(g\circ f)}="6";
(20,-30)*+{\sigma_C \circ (Fg\circ Ff)}="5";
{\ar@{=>}^{\id_{Gf}\star \sigma_f} "1"; "2"};
{\ar@{=>}^{(a')\inv} "2"; "3"};
{\ar@{=>}^{\sigma_g\star \id_{Ff}} "3"; "4"};
{\ar@{=>}_{a'} "4"; "5"};
{\ar@{=>}^{a'} "8"; "1"};
{\ar@{=>}_{\sigma_{gf}} "7"; "6"}; 
{\ar@{=>}_{\tau\star  \id_{\sigma_A}}  "8";"7"} ;
{\ar@{=>}_{\id_{\sigma_C}\star \mu}  "5";"6"} ;
\endxy 
\end{equation}
and
\begin{equation}\label{eq:App2}
\xy
(-30,10)*+{1_{GA}\circ \sigma_A}="1";
(0,10)*+{\sigma_A }="2";
(30, 10)*+{\sigma_A \circ 1_{FA}}="3";
(-30, -10)*+{ G(1_A)\circ \sigma_A}="4";
(30, -10)*+{\sigma_A\circ F(1_A) \,.}="5";
{\ar@{=>}^{l'} "1"; "2"};
{\ar@{=>}^{(r')\inv } "2"; "3"};
{\ar@{=>}^{\id{\sigma_A} \star \mu} "3"; "5"};
{\ar@{=>}_{\sigma_{1_A}} "4"; "5"}; 
{\ar@{=>}_{\tau\star \id_{\sigma_A}} "1";  "4"};
\endxy 
\end{equation}
Here $a', r'$ and $l'$ are the associator and the unitors in the
target bicategory $\calB'$.
\end{definition}

\begin{definition}
Let $F,G:\calB\to \calB'$ be two (pseud-) functors between two (2,1)
bicategories and $\sigma,\sigma': F\Rightarrow G$ two (psudo-) natural
transformations.  A {\sf modification} $\Gamma:\sigma \rightsquigarrow
\sigma'$  consists of the following data:
\begin{itemize}
\item for every object $A$ of $\calB$ a 2-cell $\Gamma_A:
  \sigma_A\Rightarrow \sigma'_A$ between 1-cells $\sigma_A,\sigma'_A: FA\to GA$.
  \end{itemize}
The data are subject to the following conditions: for any 1-cell
$f:A\to B$ of $\calB$  the diagram
\[
\xy
(-20,10)*+{G(f)\circ \sigma_A}="1";
(10,10)*+{G(f)\circ \sigma'_A }="2";
(-20, -10)*+{\sigma_B \circ F(f)}="3";
(10, -10)*+{ \sigma_B' \circ F(f)}="4";
{\ar@{<=}^{1\star \Gamma_A} "1"; "2"};
{\ar@{=>}_{\sigma_f } "1"; "3"};
{\ar@{=>}^{\sigma'_f} "2"; "4"};
{\ar@{<=}_{\Gamma_B \star 1} "3"; "4"};
\endxy 
\]
commutes.
\end{definition}

\subsection*{Functor bicategories}\mbox{}\\
Given a pair of (2,1) bicategories $\calB, \calB'$ one can define the
functor bicategory $\Hom(\calB,\calB')$ whose 0-cells are
psedo-functors, 1-cells are pseudo-natural transformations and 2-cells
are modifications; see \cite{Street}.


\begin{thebibliography}{WWWW}


\bibitem{BC} J.\ Baez and A.\ Crans, Higher-dimensional algebra VI: Lie
  2-algebras, {\em Theory and Applications of Categories} {\bf 12}
  (2004) 492--538, {\tt arXiv:math.QA/0307263.}

\bibitem{BehXu} K.\ Behrend and P.\ Xu, Differentiable stacks and
  gerbes, {\em Journal of Symplectic Geometry},{\bf 9}(3) (2011),
  285--341.

\bibitem{Benabou} J.\ B\'{e}nabou, Introduction to bicategories, in
  {\em Reports of the Midwest Category Seminar} pp. 1--77, Springer,
  Berlin, 1967.

\bibitem{Blo} C.\ Blohmann, Stacky Lie groups, {\em Int.\ Math.\ Res.\
    Notices}, 2008,
  rnn082.

\bibitem{Borceux} F.\ Borceux {\em Handbook of categorical algebra},
  volume 50 of Encyclopedia of Mathematics and its Applications,
  Cambridge Univ.\ Press, 1994.

\bibitem{SGA4} P.\ Deligne, La formule de dualit\'e globale in {\em
    Th\'eorie des topos et cohomologie \'etale des sch\'emas} Tome 3,
  Lecture Notes in Mathematics, Vol. 305. Springer-Verlag, Berlin-New
  York, 1973. vi+640 pp.

\bibitem{FW} Y.\ Fregier and F.\ Wagemann, {\em On Hopf 2-algebras},
  IMRN, vol. 2011, no. 15, 3471--3501

\bibitem{Hepworth} R.\ Hepworth, Vector fields and flows on
  differentiable stacks, {\em Theory Appl.\ Categ.} {\bf 22} (2009),
  542--587.
  
%

\bibitem{Leinster} T.\ Leinster, Basic bicategories, \url{https://arxiv.org/abs/math/9810017v1}

\bibitem{L} E.\ Lerman, Orbifolds as Stacks?. {\em L'Enseignment
Math\'{e}matique} {\bf 56} (2010), 315--363.


\bibitem{MacK} K.C.H.\ Mackenzie, {\em General Theory of Lie Groupoids
    and Lie Algebroids}, London Mathematical Society Lecture Note
  Series {\bf 213}, Cambridge University Press, Cambridge,
  2005. xxxviii+501 pp.

\bibitem{MackXu} K.C.H.\ Mackenzie and P.\ Xu, Classical lifting properties and multiplicative vector fields,
{\em Quart.\ J.\ Math.\ Oxford(2), }{ \bf 49} (1998), 59--85.

\bibitem{Metzler} D.\ Metzler, Topological and smooth stacks,
  \href{http://arxiv.org/abs/math/0306176}{arXiv:math/0306176
    [math.DG]} (2003).

\bibitem{MM} I.\ Moerdijk and J.\ Mrcun, {\em Introduction to
    Foliations and Lie groupoids}, Cambridge University Press,
  Cambridge, 2003.  ix+173 pp.


\bibitem{Noohi} B.\ Noohi, Integrating morphisms of Lie 2-algebras,
  {\em Compositio Math.} {\bf 149} (2013), 264--294.
  \href{http:doi:10.1112/S0010437X1200067Xdoi}[10.1112/S0010437X1200067X]


\bibitem{OW} C. Ortiz and J. Waldron, Lie 2-algebra of vector fields
  on Lie groupoids, {\em Journal of Geometry and Physics} {\bf  145} (2019): 103474.
 \href{http://arxiv.org/abs/1703.09791}{arXiv:1703.09791
    [math.DG]} (2017).


\bibitem{Pronk} D.A.\ Pronk, Etendues and stacks as bicategories of
  fractions, {\em Compositio Mathematica} {\bf 102}, no. 3 (1996),
  pp. 243--303.


\bibitem{SP} C. Schommer-Pries, Central extensions of smooth 2-groups and a finite-dimensional string 2-group, {\em Geometry and Topology}, {\bf 15} (2011), pp. 609-676.

\bibitem{Street} R.\ Street, Categorical Structures,  in {\em Handbook of algebra}  Vol. 1,
    pp. 529 -- 577, North-Holland, Amsterdam, 1996.
  
\bibitem{Wei} A.\ Weinstein, The Volume of a Differentiable Stack,
  {\em Lett.\ Math.\ Phys.} {\bf 90} (2009), 353--371.
\end{thebibliography}
\end{document}